\documentclass{article}
\usepackage{amsthm}

\newcounter{TheoremS} \Alph{TheoremS}
\input{amssym.tex} 

\newtheorem*{TheoremA}{Theorem A}
\newtheorem*{TheoremB}{Theorem B}
\newtheorem{Theorem}{Theorem}[section]

\newtheorem{Proposition}[Theorem]{Proposition}
\newtheorem{Lemma}[Theorem]{Lemma}

\newtheorem{Remark}[Theorem]{Remark}
\newtheorem{Assumption}{Assumption}
\newtheorem{Example}{Example}

\def\proofa{\noindent {\bf Proof of Theorem A:}}
\def\proofba{{\bf  Proof of Theorem B, Part I):}}

\def\diam{\mathop{\hbox{\rm diam}}}
\def\intset{\mathop{\hbox{\rm int}}}

\def\det{\mathop{\hbox{\rm det}}}
\def\id{\mathop{\hbox{\rm id}}}

\def\supp{\mathop{\hbox{\rm supp}}}
\def\Esup{\mathop{\hbox{\rm Esup}}}
\def\Einf{\mathop{\hbox{\rm Einf}}}
\def\osc{\mathop{\hbox{\rm osc}}}
\def\disp{\displaystyle}
\def\<<{\prec}
\def\lee{\<<}
\def\bbar{}

\def\a{\alpha}
\def\b{\beta}
\def\c{\gamma}   \def\C{\Gamma}
\def\d{\delta}
\def\D{\Delta}
\def\e{\varepsilon}
\def\k{\kappa}
\def\l{\lambda}
   
\def\s{\sigma}   
\def\t{\theta}

\def\Bbba{\chi}

\def\dd{\zeta} \def\m{\nu} \def\mathbb{\Bbb}

\def\proof{{\bf {\medskip}{\noindent}Proof: }}

\def\qed{\hfill$\square$  \bigskip} 

\begin{document}
\bibliographystyle{plain}

\title {Absolutely Continuous Invariant Measures for
Nonuniformly Expanding Maps}
\author{Huyi Hu \thanks{Mathematics Department, Michigan State University,
East Lansing, MI 48824, USA.
e-mail: $<$hu@math.msu.edu$>$.} \and Sandro Vaienti
\thanks{UMR-6207 Centre de
Physique Th\'eorique, CNRS, Universit\'es d'Aix-Marseille I,
II,Universit\'e du Sud Toulon-Var and FRUMAM, F\'ed\'ederation de
Recherche des Uniti\'es de Math\'ematiques de Marseille; address:
CPT, Luminy Case 907, F-13288 Marseille Cedex 9
e-mail: $<$vaienti@cpt.univ-mrs.fr$>$.}}

\maketitle

\begin{abstract}
For a large class of nonuniformly expanding maps of $\Bbb R^m$,
with indifferent fixed points and unbounded distorsion and  non
necessarily Markovian, we construct an absolutely continuous
invariant measure. We extend to our case  techniques previously
used for expanding maps on quasi-H\"older spaces. We give general
conditions and provide examples to which apply our result.
\end{abstract}

\section{Introduction}
\setcounter{equation}{0}

A challenge problem in smooth ergodic theory is to construct
invariant measures for multidimensional maps $T$  with some sort of
weak hyperbolicity and then to study their statistical properties
(decay of correlations, central limit theorem, distribution of
return times, etc.). For
nonuniformly expanding endomorphisms of $\Bbb R^m $, only few results
exist at the moment.  When the
system has a Bernoulli structure, or verifies the so-called
``finite range structure'', and it enjoys a suitable distortion
relation (Renyi's condition),  M. Yuri \cite{Yu1,Yu2} was able to
construct an invariant, possibly $\sigma$-finite, measure
absolutely continuous with respect to the Lebesgue measure.
Young's tower (\cite{Yo1, Yo2}), which is mainly for
nonuniformly hyperbolic systems,  also works for  nonuniformly
expanding maps, and the invariant measures and other statistical
properties can be obtained, if some bounded distortion properties
are assumed.
In Alves-Bonatti-Viana's work (\cite{ABV}, also see \cite{A1,A2,ALP}),
nonuniformly expansion are understood as the average value of
$\log ||DT(x)^{-1}||$ along the orbits to be less than zero for almost
every points.  Under some conditions on the set of critical points,
they can construct an absolutely continuous invariant measure.
Recently this theory
has been applied to maps which allow contraction in some regions
\cite{O, AT}.

The aim of our paper is to treat a class of nonsingular
transformations with indifferent fixed points which do not enjoy
any Markov property. We obtain existence of absolutely continuous
invariant measures that can be finite or infinite, depending on
the behaviour of $T$ near the fixed point. The technique we use
consists of the following steps.  We first replace the
transformation with the first return map with respect to the
domain outside a small region around the indifferent fixed point.
What we get is a uniformly expanding map with a countable number
of discontinuity surfaces.  Then we  prove a Lasota-Yorke
\cite{LY}  inequality on the induced space   by acting the
Perron-Frobenius  operator on the space of ``quasi-H\"older''
functions, particularly adapted when the invariant densities are
discontinuous. As soon as the Lasota-Yorke inequality has been
proved, simple compactness argument will allow us to apply the
Ionescu-Tulcea and Marinescu theorem to conclude that there exists
an absolutely continuous invariant measure. The space of
quasi-H\"older functions, introduced by Keller \cite{K},
developed by Blank \cite{BB} and successfully applied by Saussol
\cite{S} and successively by Buzzi \cite{B} (see also \cite{BK})
and Tsujii \cite{T} to the multidimensional expanding case,
reveals to be very useful to control the oscillations of a
function under the iteration of the PF operator across the
discontinuities of the map. The use of the more standard space of
bounded variation functions allowed as well to get absolute
invariant measures for a wide  class of piecewise expanding maps,
see, for instance \cite{GB,PGB,KZ,WC}.

In adapting to our situation the Saussol's strategy to prove the
Lasota-Yorke inequality, the difficult part comes from the
indifferent fixed points. Unlike in one dimensional case, the
maps in higher dimensional space have unbounded distortion {\it
away from} the indifferent fixed points, that is, there are
uncountably many points $x$, whose neighborhoods contain points
$y$, arbitrary close to $x$, such that the distortion of $|\det
DT|$ is unbounded along the backward orbits towards the
indifferent fixed point (see Example 1 in Section 2). This forced
us to a certain number of assumptions which basically reduce to
insure sufficiently good expanding rates in a small neighborhood
of the neutral point, and insure bounded distortion along the
curves close to radial directions (Assumption 4(b) and (c)). A
careful view at the proofs will reveal that such assumptions are
unavoidable, unless to modify deeply all the structure of the
approach. We nevertheless point out that our hypothesis could be
easily verified on some simple cases if the local behavior of the
map $T$ near the indifferent fixed points is understood. On the
other hand, it seems that other known techniques are difficult to
apply. Since we know that distortions are unbounded for the maps
we are interested in, Young's results cannot be applied
directly.  Also, the condition
$\limsup_{n\rightarrow\infty}\frac{1}{n}
\sum_{i=0}^{n-1}\log||DT(T^i(x))^{-1}||<0$ in \cite{ABV} cannot
be obtained in our case (and in fact it fails if $T$ admit a
$\s$-finite absolutely continuous invariant measures). If we
study the first return maps $\hat T$ instead, then $\|D\hat
T_x(v)\|$ can be arbitrary large for $x$ close to the
discontinuity set, and therefore the assumptions on the critical
set in \cite{ABV} are not satisfied.

The plan of the paper is the following: in Section 1 we state the
assumptions and the main theorems, A and B. Section 2 is devoted to
examples. The proofs of the main results are in Section 3 through 6.


\section{Assumptions and statements of results}
\setcounter{equation}{0}

Let $M\subset \Bbb R^m$ be a compact subset with
$\overline{\intset M}=M$ and $d$ be the Euclidean distance.
Let $\nu$ be the Lebesgue measure on $M$.  We assume $\nu M=1$.

For $A\subset M$ and $\e>0$,
denote $B_{\e}(A)=\{x\in \Bbb R^m; d(x, A)\le \e \}$.

Let $T: M\to M$ be an almost
expanding piecewise smooth map with an indifferent fixed point $p$.

We assume that $T$ satisfies the following assumptions.

\begin{Assumption} {\rm (Piecewise smoothness)}
There are finitely many disjoint open sets $U_1, \cdots, U_K$ with
$M=\bigcup_{i=1}^K \overline{U}_i$ such that for each $i$,

\begin{enumerate}
\item[{\rm (a)}] $T_i:=T|_{U_i}: U_i\to M$ is $C^{1+\a}$;

\item[{\rm (b)}] $T_i$ can be extended to a $C^{1+\a}$ map
$T_i: \tilde {U_i}\to M$ such that $T_i\tilde U_i\supset
B_{\e_1}(T_iU_i)$ for some $\e_1>0$, where $\tilde U_i$ is a
neighborhood of $U_i$.

\end{enumerate}
\end{Assumption}


\begin{Assumption} {\rm (Fixed point)}
There is a point $p\in U_1$ such that:

\begin{enumerate}
\item[{\rm (a)}] $Tp=p$;

\item[{\rm (b)}] $T^{-1}p\notin \partial U_j$ for any $j$.
\end{enumerate}
\end{Assumption}

Since $M\subset \Bbb R^m$, we may take a coordinate system such that $p=0$.
Hence, we write $|x|=d(x,p)$ if $x\in M$.

For any $x\in U_i$, we define $s(x)=s(x, T)$ by
$$
s(x, T)=\min\bigl\{s: d(x,y)\le sd(Tx, Ty),
 y\in U_i, d(x,y)\le \min\{\e_1, 0.1|x|\} \ \bigr\}.
$$

Denote by $\c_m$ the volume of the unit ball in $\Bbb R^m$.


\begin{Assumption} {\rm (Expanding Rates)}
There exists an open region $R$ bounded by a smooth surface
with $p\in R$, $\overline R\subset TR$, $\overline{TR}\subset U_1$
and with either $\overline R\subset TU_j$
or $\overline R\cap TU_j=\emptyset$ such that:

\begin{enumerate}
\item[{\rm (a)}] $0<s(x)\le 1$ \ $\forall x\in M\setminus\{p\}$,
and if $s(x)=1$ then $x\in R$ and $|Tx|>|x|$;

\item[{\rm (b)}] there exist constants
$\eta_0\in (0,1)$, $\e_2>0$
such that
$$
s^\a+ \lambda \le \eta_0 < 1,
$$
where
$$
s:=\max\{s(x): x\in M\backslash R\},
$$
\begin{eqnarray}\label{lambda}
\displaystyle
\lambda=\max\Bigl\{2\sup_{\e_0\le \e_2} \sup_{\e\le \e_0}
{G_U(\e,\e_0)\over \e^\a}\e_0^\a, \ {3s\c_{m-1}\over (1-s)\c_m} \Bigr\},
\end{eqnarray}
\begin{eqnarray}\label{GU}
\displaystyle G_U(\e,\e_0)=\sup_{x\in M}G_U(x, \e,\e_0),
\end{eqnarray}
and
\begin{eqnarray*}
\displaystyle G_U(x, \e,\e_0)=\sum_{j=1}^K
 {\m(T_{j}^{-1}B_\e(\partial TU_j) \cap B_{(1-s)\e_0}(x))
 \over \m(B_{(1-s)\e_0}(x))};
\end{eqnarray*}

\item[{\rm (c)}]
there exists $N=N_s>0$ and $\e_3>0$ such that for all
$x\in R_{\e_3}(TR\backslash R)$,
$$
s(T_1^{-N}(x), T_1^{N})
\le {s\over 5m} \Bigl({\l(1-s)^m \over 2C_\xi I^2}\Bigr)^{1/\a}
$$
for $\lambda$ given by (\ref{lambda}) and $I$ and $C_\xi$
given by Assumption 4(c).

\end{enumerate}
\end{Assumption}

\begin{Remark}\label{Rmk3a}
By Assumption 3(a), the map $T_j: U_j\to T_j(U_j)$ is
noncontracting for each $j$, and therefore it is a local
diffeomorphisms. Also, by the assumption, for any $x\in U_1$,
$T^{-n}_1x\to p$, because the set of limit points of
$\{T^{-n}_1x\}$ cannot contain any other point but $p$.
\end{Remark}

\begin{Remark}\label{Rmk3b}
Assumption 3(b) is the main assumption that requires uniformly expanding
outside $R$ and gives condition on the relations between expanding
rates and discontinuity.  We refer to \cite{S} for more details
about the meaning of $G_U(\e, \e_0)$.  (In fact, for small $\e_0$,
$G_U(\e, \e_0)\e_0/\e$ is greater than $4s\c_{m-1}/(1-s)\c_m^{-1}$
if there are at least two surfaces $\partial U_i$ meet at some point.
See Lemma 2.1 in \cite{S}.)
\end{Remark}

\begin{Remark}\label{Rmk3ba}
Assumption 3(b) implies $\m(\partial U_j)=0$ for any $j=1,\cdots,K$.
\footnote{In fact, if $\m(\partial U_j)>0$ for some $j$, then 
we take the set of the {\it density points}
$$
\D=\left\{x\in M: \lim_{\e\to 0}
\frac {\m(B_\e(x)\cap \partial U_j)}{\m B_\e(x)}=1 \ \right\}.
$$
By the Lebesque-Vitali Theorem (see, e.g. \cite{SG}, Chapter 10), 
$\m \D=\m(\partial U_j)>0$.  In particular, $\D\not=\emptyset$.
Therefore for any $x\in \D$, if $\e_0$ is sufficiently small
and $\e=(1-s)\e_0$, then 
$$
G_U(x, \e,\e_0)
\ge  {\m(T_{j}^{-1}B_\e(\partial TU_j) \cap B_{(1-s)\e_0}(x))
 \over \m(B_{(1-s)\e_0}(x))}
\ge  {\m(\partial U_j \cap B_\e(x))
 \over \m(B_\e(x))} 
$$
is sufficiently close to $1$, which contradicts to the assumtion.}
\end{Remark}

\begin{Remark}\label{Rmk3c}
We allow that $s(x, T)=1$ for some $x$ other than $p$. However we
still need some expanding rate inside $R$. This is given by
Assumption 3(c). If $s(T_1^{-N}(x), T_1^{N})$ can be arbitrarily
small by taking $N$ sufficiently large, then Assumption 3(c) is
always true.
\end{Remark}

Denote $R_0=TR\backslash R$.  Clearly, $R_0\subset U_1$
because of the choice of $R$.

\begin{Assumption} {\rm (Distortions)}
\begin{enumerate}
\item[{\rm (a)}]
There exists $c>0$ such that for any $x, y\in TU_j$ with
$d(x,y)\le \e_1$,
$$
\bigl|\det DT_j^{-1}(x)-\det DT_j^{-1}(y)\bigr|
\le c |\det DT_j^{-1}(x)|d(x,y)^\a,
$$
where $\e_1$ is given by Assumption 1(b);

\item[{\rm (b)}]
For any $b>0$, there exist $J>0$, $\e_4>0$ such that for any
$\e\in (0, \e_4]$, we can find
$0< N=N(\e)\le \infty$ with
$$
{|\det DT_1^{-n}(y)|\over |\det DT_1^{-n}(x)|}\le 1+J\e^\a  \quad
\forall y\in B_\e(x), \ x\in B_{\e_4}(R_0), \ n\in(0, N_{}],
$$
and
$$
\sum_{n=N_{}}^\infty\sup_{y\in B_\e(x)}|\det DT_1 ^{-n}(y)|
\le b\e^{m+\a}  \quad \forall x\in B_{\e_4}(R_0);
$$

\item[{\rm (c)}]
There exist constants $I>1$, $C_\xi>0$, $\e_5>0$ such that for
any $0< \e_0\le \e_5$, $n>0$, there is a finite or countable
partition $\xi=\xi_n$ of $B_{\e_0}(R_0)$  such that $\forall A\in
\xi$,  $0<\e\le \e_0$, $\diam (A\cap B_{\e_0}(\partial R_0))\le 5
m\e_0$,
\begin{equation}\label{vol}
{\nu\bigl(B_{\e}(\partial R_0)\cap A )\bigr)
 \over \nu\bigl(B_{\e_0}(\partial R_0)\cap A \bigr)}
\le C_\xi \left({\e\over \e_0}\right)^\a,
\end{equation}
whenever
$\nu\bigl(T_1^{-n}(B_{\e_0}(\partial R_0))\cap A \bigr)\not=0$, and
for any $x,y\in A$,
\begin{equation}\label{dist}
{|\det DT_1^{-n}(y)|\over |\det DT_1^{-n}(x)|}\le I.
\end{equation}
\end{enumerate}

\end{Assumption}

\begin{Remark}\label{Rmk4a}
In fact, Assumption 4(a) is a consequence of Assumption 1. 
\footnote{we note that by Assumption 1(b), the map $x\to |\det DT(x)|$
is continuous on $\overline{U}_i$ for each $i$.  Since 
$\overline{U}_i$ is compact, $|\det DT(x)|$ is bounded.
Hence, Assumption 4(a) follows from the fact
that $T$ is piecewise $C^{1+\a}$, Assumption 1(a)}
However, we state it here independently 
due to its importance for our arguments.
\end{Remark}

\begin{Remark}\label{Rmk4b}
If $T_1^{-1}$ has bounded distortion in $B_{\e_5}(R_0)$ in the sense
that for any $J_0>1$, there is $\e>0$ such that for any
$x,y\in B_{\e_5}(R_0)$ with $d(x,y)\le \e$ and for any $n>0$,
$\disp {|\det DT_1^{-n}(y)|\over |\det DT_1^{-n}(x)|}\le J_0d(x,y)^\a$,
then Assumption 4(b) and (c) are true with $\e_4=\e_5=\e_0$.
\end{Remark}

\begin{Remark}\label{Rmk4c}
Actually, by our proof  the condition
$\diam (A\cap B_{\e_0}(\partial R_0))\le 5 m\e_0$
in Assumption 4(c) can be replace by
\begin{eqnarray*}
\disp \diam T^{-n}_1(A\cap B_{\e_0}(\partial R_0))
\le s \Bigl({\l (1-s)^m\over2C_\xi I^2} \Bigr)^{1/\a}
\end{eqnarray*}
for all $n\ge N_s$, where $s$ and $N_s$ are given by
Assumption 3(b) and (c) respectively (see (\ref{diamAijk})).
\end{Remark}

\begin{Remark}\label{Rmk4d}
When we iterate the system, oscillations of the test functions are
produced by both discontinuities $\partial U_j$ and distortion of
$|\det DT|$. It is very common for an expanding system in
multidimensional space with an indifferent fixed point to have
unbounded distortion near the fixed point.  (See Example 1 in
Section~\ref{examples}). Assumption 4(b) requires that either the
distortion of $|\det DT_1^{-n}(x)|$ or 
$|\det DT_1^{-n}(x)|$ itself is small. On the other hand, if the
distortion of $|\det DT_1^{-n}(x)|$ is bounded along the radial
direction, then Assumption 4(c) holds.
\end{Remark}

\begin{TheoremA}
Suppose $T: M\to M$ satisfies Assumption 1-4. Then $T$ admits an
absolutely continuous invariant measure $\mu$ with at most
finitely many ergodic components $\mu_1, \cdots, \mu_s$ that are
either finite or $\sigma$-finite, and the density functions of
$\mu_i$ are bounded on any compact set away from $p$.   Hence,
\begin{enumerate}
\item[{$\cdot$}]
$\mu$ is finite if
$\disp \sum_{n=1}^\infty \nu(T_1^{-n} R)<\infty$.
\end{enumerate}
Moreover, if $|\det DT|$ is bounded and for any ball $B_\e(x)$ in $M$,
there exists ${\tilde N}={\tilde N}(x,\e)>0$ such that
$T^{{\tilde N}}B_\e(x)\supset M$,
then 
the density function is bounded below
by a positive number.  Hence
\begin{enumerate}
\item[{$\cdot$}]
$\mu$ is $\sigma$-finite if
$\disp \sum_{n=1}^\infty \nu(T_1^{-n} R)=\infty$.
\end{enumerate}
\end{TheoremA}

\begin{Remark}\label{RmkThmA}
We will give an example in Section 2 showing  that it is possible
for $\mu$ to have both finite and $\sigma$-finite ergodic components
simultaneously, and both contain the same indifferent fixed point $p$
in their supports.
\end{Remark}

Since Assumption 4(b) and 4(c) are difficult to
verify, we give some sufficient conditions in the next theorem.

One of the interesting cases we would discuss is the following:
there are constants $\c'>\c>0$, $C_i, C_i'>0$, $i=0,1,2$, such
that
\begin{eqnarray}
    |x|\bigl(1-C_0'|x|^\c+O(|x|^{\c'})\bigr)
\le  \!\!\!\!\!\!\!\!\! &&|T_1^{-1}x| \le
    |x|\bigl(1-C_0|x|^\c+O(|x|^{\c'})\bigr),  \label{B.0} \\
1-C_1' |x|^\c \le  \!\!\!\!\!\!\!\!\!\!&&\|DT_1^{-1}(x)\|
    \le  1-C_1 |x|^\c,  \label{B.1} \\
 C_2' |x|^{\c-1} \le \!\!\!\!\!\!\!\!\!\! &&\|D^2T_1^{-1}(x)\|
 \le  C_2  |x|^{\c-1}.  \label{B.2}
\end{eqnarray}
If $T$ satisfies all of the inequalities, then $\|DT_p\|=1$.
So $DT_p$ is either the identity or a rotation.
If $T$ satisfies the second inequalities in (\ref{B.0})-(\ref{B.2}),
then $\|DT_p\|$ may have eigenvalues greater than $1$.

In the theorem below,
we denote by $E(v_1, \cdots, v_k)$ the subspace spanned by
vectors $v_1, \cdots, v_k$, and by $E_x(S)$ the tangent space
of a submanifold $S$ at a point $x\in S$.
Also, we may use a coordinate system $(t, \phi)$ near $p$ where $t=|x|$
and $\phi\in {\Bbb S}^{m-1}$, the $m-1$ dimensional sphere.

\begin{TheoremB}
Suppose $T: M\to M$ satisfies Assumption 1-3 and 4(a).
Assumption 4(b) and 4(c) are satisfied
if the conditions in Part (I) and (II) below hold respectively.
Hence, the conclusions of Theorem A hold.

\begin{enumerate}
\item[I)] One of the following conditions holds:
\begin{enumerate}
\item[\ {i)}]
There exists a constant $\kappa\in (0,1)$ such that
$|\det DT|\ge \kappa^{-1} >1$, and a constant $\hat \a>\a$
such that $T$ is $C^{1+\hat\a}$ in a neighborhood of $p$.
In this case, $\mu$ is finite if
Assumption 4(c) also holds.

\item[\ {ii)}]
There exists an open region $\tilde R\subset R$ containing $p$
with $T_1^{-L}R\subset \tilde R$ for some $L>0$, and constants
$\c'>\c>0$, $C_0, C_1, C_2>0$ such that the second inequalities in
(\ref{B.0})-(\ref{B.2}) hold; and there exist constants $\d,
\tau >0$, $C_\d, C_\tau >0$ with
\begin{eqnarray}\label{B.5}
{1\over \c(1-\a)}-\tau < {\d-1\over m+\a}
\end{eqnarray}
such that for any $x\in R_0$, $n\ge L$,
\begin{eqnarray}\label{B.6}
|\det DT_1^{-n}(x)|\le {C_\d\over n^\d},  \qquad
\|DT_1^{-n}(x)\|\le {C_\tau\over n^\tau}.
\end{eqnarray}
\end{enumerate}

\item[II)] One of the following conditions holds:
\begin{enumerate}
\item[\ {i)}]
There is a decomposition of $TR$ into finite or
countable number of cones $\{{\cal C}_i\}$
and a partial order ``$\lee$'' on each ${\cal C}_i\cap R$
such that  $\nu(R\setminus \cup_i{\cal C}_i)=0$
and $T({\cal C}_i\cap R)={\cal C}_i\cap TR$;
$x\lee Tx$ for any $x\in{\cal C}_i\cap R$ and
for any $y\in R_0$ there is $x\in \partial R$
such that $x\lee y \lee Tx$;
$x\lee y$ implies $T^{-1}_1x \lee T^{-1}_1y$ and
$|\det DT(x)|\le |\det DT(y)|$.

\item[\ {ii)}]
Suppose $T$ is $C^{1+\c}$ and satisfies (\ref{B.0})-(\ref{B.2}) near $p$.
There are two families of cones $\{{\cal C}_x\}$ and $\{{\cal C}_x'\}$,
continuous uniformly in $(t, \phi)$, where $t\ge 0$ and
$\phi={\Bbb S}^{m-1}$ with $(t, \phi)\in TR$, in the tangent
bundle over the set $TR$ such that (a) $DT_x({\cal C}_x) \subset
{\cal C}_{Tx}$ and $DT_x({\cal C}_x') \supset {\cal C}_{Tx}'$ \
$\forall x\in R$; \ (b) there exists a positive angle $\theta_0$
such that for any $x\in TR$ and $v\in {\cal C}_x$ and $v'\in
{\cal C'}_x$, the angle between these two vectors is bounded from
below by $\theta_0$; \ (c) $\exists d>0$, such that
\begin{eqnarray}\label{B.strong}
{|\det DT_x|_{E(v, v')}|\over
\|DT_x|_{E(v)}\|\cdot \|DT_x|_{E(v')}\|}
\le 1-d|x|^\c
\end{eqnarray}
for any $v, v'\in {\cal C}_x$;
and
(d) ${\cal C}_x$ contains the position vector from $p$ to $x$ for all
$x\in TR$,
${\cal C}_x'$ contains $E_x(\partial B_{\e_{}}(R_0))$
for all $x\in \partial(B_\e(R_0))$, $0<\e\le \e_5$,
and
\begin{eqnarray}\label{B.weak}
\|DT_x|_{E(\partial(T_1^{-n}R))}\|
\le {|Tx|^{1/(1-\bbar \t)}\over |x|^{1/(1-\bbar \t)}}
\quad \forall x\in \partial(T_1^{-n}R), \ n>0
\end{eqnarray}
for some $\bbar \t$ with $(1+\c)(1-\bbar \t)>1$.


\end{enumerate}
\end{enumerate}
\end{TheoremB}

\begin{Remark}\label{RmkThmB1}
The condition in Theorem B.I).i) means that $DT_p$ has at least
one eigenvalue with absolute value greater than $1$.


The condition in Theorem B.II).ii) part (c) implies that
under $DT$, vectors in the cone ${\cal C}_x$ expands faster
than that in ${\cal C}_x'$.
\end{Remark}

\begin{Remark}\label{RmkThmB3}
If we write $DT(x)=T_0(x)+T_\c(x)+T_h(x)$, where $T_0=DT_p$, $T_\c$
satisfies $T_\c(tx)=t^\c T_\c(x)$\ $\forall t>0$ and
$|T_h(x)|=O(|x|^{\c'})$, $\c'>\c$,  then the cones
$\{{\cal C}_x\}$ and $\{{\cal C}_x'\}$ are mainly determined
by $T_\c$ as $x$ near $p$.
So it is easy to get uniformity near $t=0$.
\end{Remark}


\section{Examples}\label{examples}
\setcounter{equation}{0}

In the next example we show that near an indifferent fixed point
$p$ of a map $T: \Bbb R^m\to \Bbb R^m$, distortion may be
unbounded even away from $p$ in the sense that there is a point
$z$ such that for any neighborhood $V$ of $z$, we can find $\hat
z\in V$ such that the ratio
\begin{eqnarray}\label{2.0}
|\det DT_1^{-n}(z)|/|\det DT_1^{-n}(\hat z)|
\end{eqnarray}
is unbounded as $n\to \infty$.

\begin{Example}\label{example1}
  Define $T:\Bbb R^2\to \Bbb R^2$ in such a way around $(0,0)$ it
behaves like:
\begin{eqnarray}\label{f2.1}
T(x,y)=  \bigl( x(1+x^2+y^2), \ y(1+x^2+y^2)^2 \bigr).
\end{eqnarray}
\end{Example}
It is easy to see that
\begin{eqnarray}\label{2.2}
DT(x,y)=\left( \begin{array}{ll}
   1+ 3x^2+ y^2 + O(|z|^4) &  2xy+O(|z|^4)   \\
   4xy+O(|z|^4)    & 1+ 2x^2+ 6y^2 + O(|z|^4)
   \end{array} \right),
\end{eqnarray}
and
\begin{eqnarray}\label{2.3}
\det DT(x,y)=1+5x^2+7y^2+O(|z|^4),
\end{eqnarray}
where $z=(x,y)$ and $|z|=\sqrt{x^2+y^2}$.

Note that in this example, $T$ is locally injective and $T^{-1}$ will
denote its inverse.
Take $z'=(x_0, 0)$ and denote $z_n'=T^{-n}z'$.  By
Lemma~\ref{lemma3.1} in the next section, we
have $\disp |z_n'|\sim {1\over \sqrt{2n}}$, where $a_n\sim b_n$ means
$\disp \lim_{n\to\infty} {a_n\over b_n}=1$.  Hence by (\ref{2.3})
and Lemma~\ref{lemma3.2},
$\disp |\det DT^{-n}(z')|\le {D'\over n^{5/2}}$ for some $D'>0$.
On the other hand if we take $z''=(0, y_0)$ and denote
$z_n''=T^{-n}z''$, then $\disp |z_n''|\sim {1\over\sqrt{4n}}$ and
$\disp |\det DT^{-n}(z'')|\ge {D''\over n^{7/4}}$ for some $D''>0$.
So $\disp {|\det DT^{-n}(z'')|\over |\det DT^{-n}(z')|}\to \infty$
as $n\to \infty$.

Suppose that for every $z\not=(0,0)$, there is a neighborhood $V$
such that for all $\hat z\in V$, the ratio in (\ref{2.0}) is
bounded for all $n>0$. We take a curve from $z'$ to $z''$ that
does not contain the origin. By choosing finite cover on the
curve, we know that the ratio $|\det DT^{-n}(z'')|/ |\det
DT^{-n}(z')|$ should be bounded.  This is a contradiction.  It
means that there are some points away from $(0,0)$ at which
distortion is unbounded.

\medskip
In the next two examples we show how to get Assumption 4(b)
and 4(c) by applying Theorem B.

\begin{Example}\label{example2}
Let $T:\Bbb R^3\to \Bbb R^3$ be given by
$$
T(x,y,z)\!=\bigl( x(1+x^2+y^2+z^2),  y(1+x^2+y^2+z^2)^2,
                  z(2+x^2+y^2+z^2)^3 \bigr)
$$
as $(x,y,z)$ near the origin.
\end{Example}

Note that by similar arguments as above we know that for this map the
distortion is also unbounded away from the origin.

Since $\det DT_{(0,0,0)}=2$ and $T$ is $C^\infty$ near the origin,
by Theorem B.I).i), Assumption 4(b) is satisfied.

Let ${\cal C}_i$, $i=1, \cdots, 8$, be the eight octants in $\Bbb R^3$,
and define a partial order ``$\lee$'' by letting
$w_1=(x_1,y_1,z_1)\lee w_2=(x_2, y_2, z_2)$ if
$|x_1|\le |x_2|$, $|y_1|\le|y_2|$ and $|z_1|\le |z_2|$.
Clearly all the requirements in Theorem B.II).i) are satisfied.
So we get Assumption 4(c) as well.

\medskip

\begin{Example}\label{example3}
  Let $T:\Bbb R^2\to \Bbb R^2$ be defined as in the first example.
\end{Example}

For any $z=(x,y)$, we denote $z_n=T^{-n}z$.

Note that
$$
|z|(1+|z|^2+O(|z|^4))\le |Tz|\le |z|(1+2|z|^2+O(|z|^4)),
$$
or
$$
|z_n|\bigl(1+|z_n|^2+O(|z_n|^4)\bigr) \le |z_{n-1}| \le
|z_n|\bigl(1+2|z_n|^2+O(|z_n|^4)\bigr).
$$
So by Lemma~\ref{lemma3.1}, we have
\begin{eqnarray}\label{f2.4}
{1\over \sqrt{4(n+k)}}+O(n^{-\b'})
\le |z_n|
\le {1\over \sqrt{2(n+k)}}+O(n^{-\b'}),
\end{eqnarray}
for some integer $k$, where $\b'>1/2$.

Since (\ref{2.3}) implies that $|\det DT(z)|^{-1}\le
1-5|z|^2+O(|z|^4)$, by (\ref{f2.4}) and Lemma~\ref{lemma3.2} we get
\begin{eqnarray}\label{f2.5}
|\det DT^{-n}(z)|\le Dn^{-5/2}.
\end{eqnarray}
Also by (\ref{2.2}),
\begin{eqnarray*}
DT^{-1}(x,y) =\left( \begin{array}{ll}
   1 -3x^2 -y^2 + O(r^4) &  -2xy+O(r^4)   \\
   -4xy+O(r^4)    & 1- 2x^2- 6y^2 + O(r^4)
   \end{array} \right).
\end{eqnarray*}
So $\|DT^{-1}(z)\|\le 1-|z|^2+O(|z|^4)$, hence by Lemma~\ref{lemma3.2},
\begin{eqnarray}\label{f2.7}
\|DT^{-n}(z)\|\le D'n^{-1/2}
\end{eqnarray}
for some $D'>0$.  Now by (\ref{f2.5}), (\ref{f2.7}) and (\ref{B.6}),
we know that $\d=5/2$ and $\tau=1/2$.
Since $m=2$ and $\c=2$, we have (\ref{B.5}) if $\a=1/2$.
By Theorem~B.I).ii), $T$ satisfies Assumption 4(b).

Now we check that $T$ satisfies Assumption 4(c).  It is obvious that
we can use Theorem B.II).i).
However, we use this map to show how to apply Theorem~B.II).ii).

Note that if we take two vectors $v_0=(x,y)^*$ and $v_0'=(y,-x)^*$
at the tangent plane of $z=(x,y)$, where the asterisk
denotes transpose, then by (\ref{2.2}) we have
\begin{eqnarray*}
DT_{z}(v_0)=\left( \begin{array}{ll}
   x+ 3x^3+ 3xy^2 + O(|z|^5)   \\    y+ 6x^2y+ 6y^3 + O(|z|^5)
   \end{array} \right),   \\
DT_{z}(v_0')=\left( \begin{array}{ll}
   y+ x^2y+ y^3 + O(|z|^5)   \\    -x-2x^3-2xy^2 + O(|z|^5)
   \end{array} \right).
\end{eqnarray*}
This means that $|DT_{z}(v_0')|< |DT_{z}(v_0)|$.
We define ${\cal C}_z$ at each point $z$ as the cone bounded by lines
generated by vectors $3v_0+2v_0'$ and $3v_0-2v_0'$ and containing $v_0$,
and define ${\cal C}_z'$ as the cone bounded by lines
generated by vectors $3v_0'+2v_0$ and $3v_0'-2v_0$ and disjoint with
${\cal C}_z$.
We can check that Part (a) and (b) in Theorem B.II).ii) are satisfied.
Also we can check that for all unit vector $v'\in{\cal C}_z'$,
$|DT_z(v')|\le |Tz|^{2.5}/|z|^{2.5}$.
So if we take $R$ in such a way that the tangent lines of
$\partial(T^{-n_0}_1R)$ are in the cones ${\cal C}'$ for some $n_0\ge 0$,
then we use the fact $DT^{-1}({\cal C'})\subset {\cal C'}$ to get
that Part (c) is satisfied for all $n\ge n_0$ with $1-\bbar \t=2/5$.

\medskip
In the next example the absolute continuous invariant measure
$\mu$ has a finite and a $\sigma$-finite ergodic components
simultaneously, and both contain the same indifferent fixed point $p$
in their supports.

\begin{Example}\label{example4}
Suppose the map $T: M\to M$ satisfies Assumption 1 - 4(a), and in
a neighborhood, say $B_1(p)$, of the indifferent fixed point $p$,
$T$ has the form as in (\ref{f2.1}).

We also assume that there is a partition of $M=\{M_1, M_2\}$ such
that for $i=1,2$, $TM_i=M_i$ and for any ball $B_\e(x)$ in $M_i$,
there exists an integer $N$ such that $T^N B_\e(x)=M_i$, and
$$
\{z=(x,y)\in B_1(p): y< x^2\}\subset M_1, \quad \{z=(x,y)\in
B_1(p): y> x^2\}\subset M_2. \quad
$$
This is possible since it is easy to check that
$T\Gamma\cap B_1(p)=\Gamma$,
where $\Gamma=\{(x,y)\in B_1(p): y = x^2\}$.
\end{Example}

By the above example, we know that $T$ also satisfies Assumption 4(b)
and 4(c).
Therefore Theorem A can be applied.  Since both $M_1$ and $M_2$ are
invariant sets, $T$ has absolutely continuous invariant measures
$\mu_1$ and $\mu_2$ with respect to the Lebesgue measure restricted to
$M_1$ and $M_2$ respectively.  Now we show $\mu_1 M_1<\infty$ and
$\mu_2 M_2=\infty$.

For this purpose we may assume that $R=B_1(p)$.  By (\ref{f2.4}), we
know that $T_1^{-n}R\subset B_{2/\sqrt{2n}}(p)$ for all large $n$.  So
$$
\m(T_1^{-n}R\cap M_1)
\le \m\bigl\{(x,y): x^2+y^2\le {4\over 2n}, \ |y|\le |x|^2\bigr\}
\le C \Bigl({4\over 2n}\Bigr)^{3/2}
$$
for some $C>0$.  It follows that $\disp \sum_{n=1}^\infty \nu(T_1^{-n}
R\cap M_1)<\infty$.  Applying Theorem A to the system $T: M_1\to M_1$,
we get that $\mu_1 M_1\le \infty$.

Also, by (\ref{f2.4}), we have that $T_1^{-n}R\supset
B_{1/2\sqrt{4n}}(p)$ for all large $n$.  Hence it is easy to see that
$\m(T_1^{-n}R)\ge \pi/16n$ and therefore $\disp \sum_{n=1}^\infty \m
(T_1^{-n} R) = \infty$.  Since $\m \bigl(T_1^{-n}R\cap M_1 \bigr)+ \m
\bigl(T_1^{-n}R\cap M_2 \bigr) =\m \bigl(T_1^{-n}R\bigr)$, we get
$\disp \sum_{n=1}^\infty \nu(T_1^{-n} R\cap M_2)=\infty$.
So we have $\mu_2 M_2=\infty$.

\section{Proof of Theorem B, Part I)}
\setcounter{equation}{0}

We first prove a few Lemmas.

For $\c>0$, let $\b=1/\c$.
\begin{Lemma}\label{lemma3.1}
If
\begin{eqnarray}\label{3.1}
t_{n-1}\ge t_n+Ct_n^{1+\c}+O(t_n^{1+\c'}) \qquad \ \forall n>0,
\end{eqnarray}
where $\c'>\c$,  then for all large $n$,
\begin{eqnarray}\label{3.2}
t_n\le{1\over (\c C (n+k))^\b}+O\bigl({1\over  (n+k)^{\b'}}\bigr)
\qquad \ \forall n>0
\end{eqnarray}
for some $\b'>\b$ and $k\in {\Bbb Z}$.  The result remains true if
we exchange ``$\le$'' and ``$\ge$''.  Therefore, if (\ref{3.1})
becomes an equality, then so does (\ref{3.2}).
\end{Lemma}

\proof
We claim that if
\begin{eqnarray}\label{3.1a}
t_{n-1}\ge t_n+Ct_n^{1+\c}+C't_n^{1+\c'},
\end{eqnarray}
for some large $n$ and
\begin{eqnarray}\label{3.1b}
t_n^\c \ge {1\over \c Cn}\Bigl(1+{1\over n^{\d'}}\Bigr)
\end{eqnarray}
for some $\d' >0$,  then
\begin{eqnarray*}
t_{n-1}^\c \ge {1\over \c C(n-1)}\Bigl(1+{1\over (n-1)^{\d'}}\Bigr).
\end{eqnarray*}
This gives the results since we can choose an integer $k$ such that
for some large $n_0>0$,
$$t_{n}^\c \le {1\over \c C(n_0+k)}\Bigl(1+{1\over (n_0+k)^{\d'}}\Bigr).$$
By relabelling the indices, the claim implies (\ref{3.2})
for all $n\ge n_0$.

Now we prove the claim.
Denote $\c_n=\c\bigl(1+n^{-\d'}\bigr)^{-1}$.
By (\ref{3.1a}) and (\ref{3.1b}),
$$t_{n-1}^\c
\ge {t_n}^\c\bigl(1+C{t_n}^\c+C' {t_n}^{\c'}\bigr)^\c \ge {1\over
Cn \c_n} \Bigl(1+{C\over C n \c_n}
       +{C'\over (Cn\c_n)^{\c'/\c}}\Bigr)^\c.$$
To prove the lemma we only need to show that
$${1\over n\c_n}\Bigl(1+{{1\over n\c_n} }
      +{C'\over (Cn\c_n)^{\c'/\c}}\Bigr)^\c
\ge {1\over (n-1)\c_{n-1}},$$
or, equivalently,
$${n-1\over n}\Bigl(1+{1\over n\c}+{1\over n^{1+\d'}\c}
      +{C'\over (Cn\c_n)^{\c'/\c}}\Bigr)^\c
\ge {\c_n\over \c_{n-1}}
 =  {1+(n-1)^{-\d'}\over 1+n^{-\d'}}.  $$
Take $\d'< \min\{1, \c'/\c-1\}$.  Then $(n\c_n)^{-(\c'/\c)}$
is of higher order.
We can check that as $n\to \infty$, the left side of the inequality
is like $1+n^{-(1+\d')}$
and the right side is like $1+\d' n^{-(1+\d')}$.
Since $\d'< 1$, the right side is smaller as $n$ large.
\qed

\begin{Lemma}\label{lemma3.2}
If for all $n>0$, $t_n$ satisfies (\ref{3.2}), and
$r(t_n)\le 1-C't_n^\c+O(t_n^{1+\c'})$, where $C'>0$,  then there exists
$D>0$ such that for all $k_0\ge k$,
\begin{eqnarray}\label{3.3}
\prod_{i=k_0-k}^{n+k_0-k} r(t_i)\le D \Bigl({k\over n+k}\Bigr)^{C'/\c C}.
\end{eqnarray}
The result remains true if we replace ``$\le$'' by ``$\ge$''
in all three inequalities.
\end{Lemma}

\proof
Note that
$$
r(t_n)
\le 1-{C'\over \c Cn}+O\Bigl({1\over n^{1+\c'}}\Bigr)
=  \Bigl(1-{1\over n}\Bigr)^{C'\over \c C}
  \cdot \Bigl(1+O\bigl({1\over n^{1+\c'}}\bigr)\Bigr),
$$
where $\c'>0$.   Then we take the product. \qed

\begin{Lemma}\label{lemma3.3}
Let $\t\in (0,1)$ and $\bar C_1', \bar C_2, \bar D_1 > 0$,
and let $\tilde R\subset \Bbb R^m$ be a bounded
region containing the origin.
Suppose the map $T: \tilde R\to \Bbb R^m$ is injective with
$T^{-1}\tilde R\subset \tilde R$ and satisfies
\begin{eqnarray}
d(Tx, Ty)\ge (1+\bar C_1'|x|^\c) d(x, y),  \label{3.4}\\
\log\Bigl|{\det DT(x)\over\det DT(y)}\Bigr| \label{3.5}
\le \bar C_2|x|^{\c-1}d(x,y)
\end{eqnarray}
for all $x,y\in \tilde R$ with $d(x,y)\le |x|/2$.
Then there exists ${J'}>0$ such that for all $x,y\in T\tilde R$ with
\begin{eqnarray}\label{3.6}
d(x_i, y_i)^{1-\t}\le \bar D_1 |x_i|,   \quad i=1,\cdots, n,
\end{eqnarray}
where $x_i=T^{-i}x$ and $y_i=T^{-i}y$, we have
\begin{eqnarray}\label{3.7}
\log\Bigl|{\det DT^n(x_n)\over\det DT^n(y_n)}\Bigr|
\le {J'}d(x,y)^\t.
\end{eqnarray}
\end{Lemma}

\proof
We prove by induction that for all $i=1,\cdots,n$,
\begin{eqnarray}\label{3.8}
\log\Bigl|{\det DT^{i}(x_n)\over\det DT^{i}(y_n)}\Bigr|
\le {J'}d(x_{n-i},y_{n-i})^\t.
\end{eqnarray}

For $i=1$, by (\ref{3.5}), (\ref{3.6}) and (\ref{3.4}), we have
\begin{eqnarray*}
\log\Bigl|{\det DT(x_n)\over\det DT(y_n)}\Bigr|
\le \bar C_2\bar D_1|x_n|^\c d(x_{n},y_{n})^\t
\le \bar C_2\bar D_1|x_{n-1}|^\c d(x_{n-1},y_{n-1})^\t.
\end{eqnarray*}
So if ${J'}\ge \sup\{\bar C_2\bar D_1|x|^\c: \ x\in \tilde R\}$
then the right side of the inequality is less than
${J'} d(x_{n},y_{n})^\t$ because $|x_n|\le |x|$.

Suppose (\ref{3.8}) is true up to $i=k-1$.
Then similarly we have
\begin{eqnarray*}
  & &\log\Bigl|{\det DT^{k}(x_n)\over\det DT^{k}(y_n)}\Bigr|
\le \log\Bigl|{\det DT^{k-1}(x_n)\over\det DT^{k-1}(y_n)}\Bigr|
 +  \log\Bigl|{\det DT(x_{n-k+1})\over\det DT(y_{n-k+1})}\Bigr|    \\
&\le &{J'}d(x_{n-k+1},y_{n-k+1})^\t
    + \bar C_2|x_{n-k+1}|^{\c-1}d(x_{n-k+1},y_{n-k+1}) \\
&=  &{J'}\Bigl( 1+{\bar C_2\bar D_1\over {J'}} |x_{n-k+1}|^\c\Bigr)
      \cdot {d(x_{n-k+1},y_{n-k+1})^\t \over d(x_{n-k},y_{n-k})^\t}
      \cdot d(x_{n-k},y_{n-k})^\t    \\
&\le &{J'}\Bigl( 1+{\bar C_2\bar D_1\over {J'}} |x_{n-k+1}|^\c\Bigr)
      \cdot {1\over (1+\bar C_1'|x_{n-k+1}|^\c)^\t}
        d(x_{n-k},y_{n-k})^\t.
\end{eqnarray*}
Clearly if ${J'}$ is large enough, then the right side is bounded by
${J'} d(x_{n-k},y_{n-k})^\t$.  We get (\ref{3.8}) for $i=k$.
\qed

\proofba

i) \
We may assume that $T$ is $C^{1+\hat\a}$ and
$|\det DT|\ge \kappa^{-1}>1$
on $TR$, because otherwise we can increase $N(\e)$ and $J$.
We may also regard $\hat\a\le 1$.
So there exist $c_1>0$ such that
$$
{|\det DT_1^{-1}(y)| \over |\det DT_1^{-1}(x)|}
\le 1+c_1d(x,y)^{\hat\a}
$$
for all $x, y\in TR$.  Let $x_i=T_1^{-i}x$ and
$y_i=T_1^{-i}y$.  Clearly, $d(x_i, y_i)\le d(x,y)$.
So if $d(x,y)\le \e$ and $0<n\le N$, then
\begin{eqnarray}\label{B.aa1}
{|\det DT_1^{-n}(y)| \over |\det DT_1^{-n}(x)|}
\le \bigl( 1+c_1d(x,y)^{\hat\a} \bigr)^n
\le \bigl( 1+c_1\e)^{\hat\a} \bigr)^N
\end{eqnarray}
Also, there exists $C>0$ such that for any
$y\in B_\e(R_0)$,
$|\det DT_1^{-n}(y)|\le C\kappa^n$.
Hence,
$$
\sum_{n=N}^\infty\sup_{y\in B_\e(x)}|\det DT_1 ^{-n}(y)|
\le {C\kappa^N\over 1-\kappa}.
$$

Let $b>0$ be given.

Consider the function
$$
\s(\e)={( 1+c_1\e^{\hat\a})^{N_0-c_2\log\e} \over 1+J\e^\a},
$$
where $N_0=1+\log(C^{-1}b(1-\k))/\log \k$ and
$c_2=-(m+\a)/\log \k$.
Since $\disp \lim_{\e\to 0}(1+c_1\e^{\hat\a})^{N_0-c_2\log\e}=1$,
we have $\disp \lim_{\e\to 0}\s(\e)=1$.
Note that if
\begin{eqnarray}\label{B.aa2}
  (N_0-c_2\log\e)\hat\a c_1\e^{\hat\a-1}\cdot (1+J\e^\a)
 -\a J\e^{\a-1}\cdot ( 1+c_1\e^{\hat\a}) <0,
\end{eqnarray}
then $\s'(\e)<0$.
Since $\hat\a>\a$, the first term in (\ref{B.aa2}) is of higher
order.  So we can choose $J>0$ and $\e_4>0$ such that
for all $\e\in(0, \e_4]$, (\ref{B.aa2}) holds and therefore
$\s(\e)\le 1$.

Now for each $\e\in (0,\e_4]$, we take $N=N(\e)$ as the
integer part of $N_0-c_2\log\e$.  Clearly, for such $N$ we have
$$
{C\kappa^N\over 1-\kappa}\le b\e^{m+\a}.
$$
So the second inequality in Assumption 4(b) is true.
For the first inequality, note that
$$
(1+c_1\e^{\hat\a})^{N}
\le ( 1+c_1\e^{\hat\a})^{N_0-c_2\log\e}
\le 1+J\e^{\a}.
$$
Then by (\ref{B.aa1}) we get what we need.

ii) \ Denote $\b=1/\c$ and $\t=\a$.  Take $\rho>0$ such that
\begin{eqnarray}\label{B.11}
{\b\over 1-\t}-\tau< \rho < {\d-1\over m+\a}.
\end{eqnarray}

Let $b>0$ be given.

Note that by Lemma~\ref{lemma3.1}, (\ref{B.0}) implies
that there exists $\bar C_0>0$ such that for any $x\in R_0$,
$\disp |x_n|\ge {1\over (\bar C_0 n)^\b}$.
Take $N_b\ge L$ such that for all $n\ge N_b$,
\begin{eqnarray}\label{B.13}
  b^{-{1\over m+\a}}
  \Bigl(\sum_{k=n}^\infty {C_\d\over k^\d} \Bigr)^{{1\over m+\a}}
< {1\over n^\rho}
< {1\over (n-1)^\rho}
< {1\over 2C_\tau \bar C_0^{{\b\over 1-\t}}n^{{\b\over 1-\t}-\tau}},
\end{eqnarray}
where $C_\d$ and $C_\tau$ are as in (\ref{B.6}).
The inequality is possible because of (\ref{B.11}).

Note that (\ref{B.1}) and (\ref{B.2}) imply (\ref{3.4}) and
(\ref{3.5}) respectively.
By Lemma~\ref{lemma3.3} we can take ${J'}>0$ such that
(\ref{3.7}) holds for any $x\in R_0$, $n>0$
whenever (\ref{3.6}) holds with $\bar D_1=1$ for all $x_i, y_i$,
$i=1,\cdots, n$.

Take $\e_4'>0$ such that for all $x, y$ with $x\in R_0$,
$d(x,y)\le \e_4'$, $n=1,\cdots, N_b$, we have
$d(x_n, y_n)^{1-\t}\le |x_n|$.  By the choice of ${J'}$,
(\ref{3.7}) holds for all $1\le n\le N_b$.

Then we take $\e_4=\min\{\e_4', 1/N_b^\rho\}$,
and $J>0$ such that $e^{J'\e_4^\t}\le 1+J\e_4^\t$.

We show that $J$ and $\e_4$ satisfies the requirement.
Let $\e\in (0, \e_4]$.  Take $N=N(\e)> N_b$ such that
$$
{1\over N^\rho} \le \e < {1\over (N-1)^\rho}.
$$
By the first inequality of (\ref{B.6}) and (\ref{B.13}),
$$
    \sum_{k=N}^\infty \sup_{y\in B_\e(x)} |\det DT^{-k}(y)|
\le \sum_{k=N}^\infty {C_\d\over k^\d}
\le b\cdot {1\over N^{\rho(m+\a)}}
\le b\e^{m+\a}.
$$
On the other hand, if $x\in R_0$ and $d(x,y)\le \e$,
then by the last inequality of (\ref{B.6}) and (\ref{B.13}),
for any $N_b< n\le N$,
$$
    d(x_n, y_n)
\le {2C_\tau\over n^\tau}\e
\le {2C_\tau\over n^\tau}{1\over (N-1)^\rho}
\le {1\over {\bar C_0}^{{\b\over 1-\t}}n^{{\b\over 1-\t}}}
\le |x_n|^{{1\over 1-\t}}.
$$
So we know that (\ref{3.7}) holds for all $0\le n\le N$.
Then by the choice of $J$ and the fact $\t=\a$,
$$
   \Bigl|{\det DT^{n}(x_n)\over\det DT^{n}(y_n)}\Bigr|
\le e^{{J'}d(x,y)^\t}
\le e^{{J'}\e^\a}
\le 1+J\e^\a.
$$
This is what we need.
\qed

\section{Proof of Theorem B, Part II)}
\setcounter{equation}{0}

This proof consist of two parts, i) and ii).

{\bf i)} \ For $x\in \partial R$, denote
$$
{\cal D}(x)=\{z\in R_0: \ x\lee z\lee Tx\}.
$$
Clearly the collection
$\{{\cal D}(x): \ x\in \partial R\cap{\cal C}_i \}$
form a cover of ${\cal C}_i\cap R_0$.  So we can construct
a partition $\xi$ of $R_0$ such that every element
of $\xi$ belongs to some ${\cal D}(x)$.

Note that for any $x$,
$${|\det DT_1^{-n}(x)|\over |\det DT_1^{-n}(Tx)|}
= {|\det DT(x)|\over |\det DT(x_n)|}
\le |\det DT(x)|
$$
is always bounded.  So for any $y,z\in {\cal D}(x)$,
we have
$$
{|\det DT_1^{-n}(y)|\over |\det DT_1^{-n}(z)|}
\le {|\det DT_1^{-n}(x)|\over |\det DT_1^{-n}(Tx)|}
\le |\det DT(x)|.
$$
Hence (\ref{dist}) follows.  Obviously we can arrange the
partition $\xi$ in such a way that (\ref{vol}) also holds.
Therefore $\xi$ is a desired partition for any $n$.

{\bf ii)} \ First, we take $\bbar \t>0$ such that
$$
DT_x(v')\le (|Tx|/ |x|)^{1/(1-\bbar \t)}
$$
for all $x\in \partial(T^{-n}_1R)$ and $v'\in
E_x(\partial(T^{-n}_1R))$. This is possible because of the
assumption stated in Part (d) of Theorem~B.(II).
So for any $n>0$, if we take $x,y\in \partial R_{0}$
such that $d(x_n, y_n)\le \bar D_1|x_n|^{1/(1-\t)}$, we have
\begin{eqnarray}\label{distanceb}
d(x_i, y_i)\le \bar D_1|x_i|^{1/(1-\t)} \qquad
\forall i=1,\cdots, n.
\end{eqnarray}
By Lemma~\ref{lemma3.3}, we get that there exists $I_1>0$ such that
\begin{eqnarray}\label{distortionb}
{|\det DT_1^{-n}(y)|\over |\det DT_1^{-n}(x)|}\le I_1.
\end{eqnarray}
That is, (\ref{dist}) holds for all such $x,y$.

We construct $\xi=\xi_n$.
Note that we only need do it for $n$ sufficiently large.
Since the family of cones ${\cal C}'_x$ are continuous uniformly
in $(t, \phi)$, we can find $t_0>0$ such that for any $x, y\in TR$
with $d(x,y)\le t_0$, the Hausdorff distance between ${\cal C}'_x$ and
${\cal C}'_y$ is less than $\theta_0/2$.  Then we take $N>0$
large enough such that for any $x\in R_0$ and $n>N$,
$|x_n|\le t_0$.
Note that for any $x$, the position vector from $p$ to $x$,
denoted by $u_x$, is contained in ${\cal C}_x$.
By Part (a) and (d) in the conditions of the theorem we know that
at $x\in T_1^{-n}(\partial R_0)$,
${\cal C}_{x}'$ contains the tangent plane of the surface.
Hence, if $v'\in E_x(T_1^{-n}(\partial R_0)$,
then the angle between $u_x$ and $v'$, denoted by $\angle(u_x, v')$,
is larger than $\theta_0$, and therefore for any
$v'\in E_y(T_1^{-n}(\partial R_0)$, we have $\angle(u_x, v')\ge \theta_0/2$,
whenever $y\in T_1^{-n}(\partial R_0)$ with $d(x,y)\le t_0$.
So for any $x, y\in T_1^{-n}(\partial R_0)$ with $d(x,y)\le t_0$,
we have $d_S(x,y)\le d(x,y)/\sin(\theta_0/2)$, where $d_S(\cdot,\cdot)$
is the distance restricted to the surfaces $\{T_1^{-n}(\partial R_0)\}$.
This means that we can take a partition $\xi^{(n)}$ on
$T_1^{-n}(\partial R_0)$ such that every element of $\xi^{(n)}$
is contained in a ball of radius $|x_n|^{1/(1-\t)}$
and containing a ball of radius $|x_n|^{1/(1-\t)}/10m\sin(\theta_0/2)$,
with respect to the metric on $T_1^{-n}(\partial R_0)$, and these
elements are close to $(m-1)$ dimensional disks.
Denote $\xi'=T^n \xi^{(n)}$.  Clearly, it is a partition
of $\partial R_0$.
Then we can take a partition $\xi$ of $R_0$ whose elements has
the form $\disp \cup_{x\in A'} {\cal F}_x\cap R_0$, where $A'$ is
an element of $\xi'$, and ${\cal F}_x$ is given in Lemma~\ref{lemma4.0}.

Now we prove that $\xi$ satisfies (\ref{vol}) and (\ref{dist}).
Condition (\ref{B.1}) implies $\|DT(p)\|=1$.
We first consider the case that $DT(p)=\id$.

By (\ref{B.0}), we know that $d(x, Tx)\le C|x|^{1+\c}$ for some
$C>0$.  So the ``width'' of the annulus
$T_1^{-i}(B_{\e_{5}}(R_0))$ is bounded by $C'|T_1^{-i}x|^{1+\c}$
for some $C'>0$. By Part (b) and (d) of the condition in the
theorem, for $0<\e\le \e_5$, $x\in T_1^{-i}(\partial B_\e(R_0))$,
the angle between the tangent space of $T_1^{-i}(\partial
B_\e(R))$ and the position vector $u_x$ is larger than
$\theta_0$. So the length of the curve ${\cal F}_{T_1^{-i}x}\cap
T_1^{-i}B_{\e_{5}}(R_0)$ is bounded by $C|T_1^{-i}x|^{1+\c}$ for
some $C\ge C'$. Hence, for any $x, y\in B_{\e_{5}}(R_0)$ with
$y\in {\cal F}_x$, we can get
\begin{eqnarray}\label{distancec}
d(x_i, y_i)\le C|x_i|^{1+\c}
\end{eqnarray}
and therefore by applying Lemma~\ref{lemma3.3} get
\begin{eqnarray}\label{distortionc}
{|\det DT_1^{-n}(y)|\over |\det DT_1^{-n}(x)|}\le I_2
\end{eqnarray}
for some $I_2>0$.
Also, the construction of $\xi'$ implies (\ref{distanceb})
and therefore (\ref{distortionb}) for any $x,y\in A'$,
where $A'\in \xi'$.  So by the construction of
$\xi$, we get (\ref{dist}) with $I=I_1I_2^2$ for any
$x,y\in A$.

On the other hand, for any $x, y\in B_{\e_{5}}(R_0)$ with $y\in
{\cal F}_x$, we have (\ref{distancec}).  So we can apply
Lemma~\ref{lemma4.3} to get that inside $A$, distortion of
$|DT|_{E({\cal F})}|$ is bounded. It means that for each $x\in A$,
the ratio of the length of $T_1^{-n}(B_{\e}(\partial R_0)\cap A
))\cap {\cal F}_{x_n}$ and the length of
$T_1^{-n}(B_{\e_0}(\partial R_0)\cap A ))\cap {\cal F}_{x_n}$ is
uniformly bounded by $\e/\e_0$ multiplied by a constant. Notice
that the angle between the tangent vectors of ${\cal F}$ and the
tangent space of $T_1^{-n}(\partial B_{\e}(R_0))$ are greater
than $\theta_0$.  Also notice that by the construction of $\xi'$,
for any $A\in \xi$, the size of the set $T_1^{-n}A$ along the
fiber direction is much smaller than the size of $T_1^{-n}A'$.
Hence, the ratio between $\nu(T_1^{-n}B_{\e}(\partial R_0)\cap A
))$ and $\nu(T_1^{-n}B_{\e_0}(\partial R_0)\cap A ))$ is bounded
by a constant times $\e/\e_0\le (\e/\e_0)^\a$ for some $\a\in
(0,1]$. Now we use (\ref{dist}) to get (\ref{vol}).\footnote{Let
us make this argument more precise. We denote with $A'_n$ and
$A_n(\e)$ respectively the backward iterates  $T_1^{-n}A'$ of
some $A'\in \xi'$ and of the set $A\cap B_{\e}(\partial R_0)$
where $A = \disp \cup_{x\in A'} {\cal F}_x\cap R_0$. Since the
angles between the tangent spaces of the curves ${\cal F}_x$ and
the tangent spaces of the $\epsilon$-neighborhood of the boundary
of $R_0$ are uniformly bounded away from zero, the length of the
curve ${\cal F}_x \cap B_{\e}(\partial R_0)$, when $x\in A'$, is
of order $\epsilon$. Its $n$-backward iterate in $A_n(\e)$ will
be therefore bounded by a constant times $\epsilon$ times
$d_{n,M}^{1+\c}$, where $d_{n,M}$ is the maximum over the
$\epsilon$-compact neighborhood of $R_0$ of $|T_1^{-i}x|$ (see
above; equivalently we set $d_{n,m}$  the minimum of
$|T_1^{-i}x|$ over the $\epsilon$-compact neighborhood of $R_0$).
Let us call this upper bound $l_{n,\epsilon}$. We construct then
the $l_{n,\epsilon}$-neighborhood
 of $A'_n$, $B_{l_{n,\epsilon}}(A'_n)$. Clearly
 $$
\nu(A_n(\e))\le \nu(B_{l_{n,\epsilon}}(A'_n))\le
\nu'(A'_{n,\e})l_{n,\epsilon}
 $$
where $A'_{n,\e}= \{z\in T_1^{-n}A'; d(z, A'_n)\le
l_{n,\epsilon}\}$ and $\nu'$ denotes the riemannian volume  on
$T_1^{-n}\partial R_0$. Since $A'_{n,\e}$ is contained in a ball
of radius $d_{n,M}^{\frac{1}{1-\theta}}+l_{n,\epsilon}$ and
$A'_n$ by construction contains a ball of radius
$\frac{d_{n,m}^{\frac{1}{1-\theta}}\sin \theta_{0}/2}{10m}$, we
have that $\nu'(A'_{n,\e})\le \mbox{const}
(d_{n,M}^{\frac{1}{1-\theta}}+l_{n,\epsilon})^{m-1} \gamma_{m-1}$
and $\nu'(A'_n)\ge (\frac{d_{n,m}^{\frac{1}{1-\theta}}\sin
\theta_{0}/2}{10 m})^{m-1}\gamma_{m-1}$. But $d_{n,M}, d_{n,m}$
are of order $n^{-\beta}$, with $\beta=1/\gamma$ (see Lemma 3.1),
and since $(1+\gamma)(1-\theta)>1$, we see immediately that for
large $n$:
$$\nu(A_n(\e))\le C' \nu'(A'_n)l_{n,\epsilon} $$ where $C'$ is a suitable
constant, depending on $m$. Let us now define the following
objects: $A_n(\e_0)$: the backward iterate of $A\cap
B_{\e_0}(\partial R_0)$, $l'_{n,\epsilon_{0}}$: the minimum
length of the backward images of the curves ${\cal F}_x \cap
B_{\e_0}(\partial R_0)$, when $x\in A'$; $A'_{n,\e_0}=\{z\in
A'_n; d(z,\partial A'_n)\ge l'_{n,\epsilon_{0}}\}$ and
$B_{l''_{n,\e_0}}(A'_{n,\e_0})$ the
$l''_{n,\epsilon_{0}}$-neighborhood
 of $A'_{n,\e_0}$, being $l''_{n,\epsilon_{0}}=l'_{n,\epsilon_{0}}\sin \theta_0$. Moreover by what we already said
 above and which follows from Lemma 4.3, the bounded distortion property along the points of the
backward images of the curves $A\cap B_{\e_0}(\partial R_0)$,
will imply  that $l'_{n,\epsilon_{0}}$ will be of the same order
as $l_{n,\epsilon_{0}}$ (the maximum length of the backward images
of the curves). Taking this into account we get:
$$
\nu(A_n(\e_0))\ge \nu(B_{l''_{n,\e_0}}(A'_{n,\e_0}))\ge
((\frac{d_{n,m}^{\frac{1}{1-\theta}}\sin
\theta_{0}/2}{10m}-l''_{n,\epsilon_{0}})^{m-1} \gamma_{m-1}
l''_{n,\epsilon_{0}}
$$
By using as above the uniform bounds on $d_{n,M}, d_{n,m}$ when
$n$ is large, we see that $\nu(A_n(\e_0))\ge C''
\nu'(A'_n)l_{n,\epsilon_{0}}$, where $C''$ is a suitable constant
depending on $m$. By dividing $\nu(A_n(\e))$ and
$\nu(A_n(\e_0))$, we get the desired result.}

If $DT_p\not= \id$, then it is a rotation, say $S$.  Hence near $p$
we can write $Tx=Sx+T_r(x)$ where $|T_r(x)|\le C|x|^{1+\c}$.  If
we write $T^{(i)}=\id+S^{-i}\circ T_r \circ S^{i-1}$, then
$T^n=S^n\circ T^{(n)}\circ\cdots \circ T^{(1)}$.
It implies that   the ``width'' of the annulus $T_1^{-i}R_0$
is bounded by   $C|T_1^{-i}x|^{1+\c}$.
Then we apply the same arguments to get (\ref{dist}) and (\ref{vol}).
\qed

\begin{Lemma}\label{lemma4.0}
There is a foliation on $\{{\cal F}_x\}$ on $TR\setminus \{p\}$
consisting of curves from $p$ to points on $\partial(TR)$
such that for any $x\in TR$, the tangent line of ${\cal F}_x$
lies in ${\cal C}_x$, and $T{\cal F}_x\cap TR={\cal F}_{Tx}$.
\end{Lemma}

\proof
Denote $\disp E_x=\cap_{n\ge 0}DT_{T_1^{-n}x}^n({\cal
C}_{T_1^{-n}x})$ for all $x\in TR\setminus \{p\}$.  By
Lemma~\ref{lemma4.1}, we know that sine of the angle between any
two vectors in $DT_{T_1^{-n}x}^n({\cal C}_{T_1^{-n}x})$
is less than $(1-d|x_n|^\c)\cdots (1-d|x_1|^\c)$. By (\ref{B.0}) and
Lemma~\ref{lemma3.1}, the product diverges as $n\to \infty$. So
$\{E_x\}$ is a subbundle of the tangent bundle over $TR \setminus
\{p\}$.  Further, we have $DT_x(E_x)=E_{Tx}$ for all $x\in R$. By
Lemma~\ref{lemma4.2}, we know that $\{E_x\}$ satisfies the
H\"older condition near each $x$ with H\"older constants
depending on $x$. Note that $\{E_x\}$ determines a vector field.
We can integrate it to get a family of curves $\{{\cal F}_x\}$
from $p$ to boundary points of $TR$ that satisfies $T{\cal
F}_x\cap TR={\cal F}_{Tx}$. By our assumption, $\{{\cal F}_x\}$
is the ``strong unstable manifold'' at $x$.

It is easy to see that
the curve passing through $x$ is unique,  and therefore $\{{\cal
F}_x\}$ forms a foliation. In fact, if there are two such curves
${\cal F}_x$ and ${\cal F}_x'$ that pass through $x$, then we can take
a curve $\C$ close to $x$ joining $y\in {\cal F}_x$ and
$y'\in{\cal F}_x'$ such that the tangent line of $\C$ is in
${\cal C}'$. Let us denote by $A_n$ the area of the ``triangle''
bounded by the curves $T^{-n}_1\C$, $T^{-n}_1{\cal F}_{x,y}$ and
$T^{-n}_1{\cal F}'_{x,y'}$, and by $L_n$ and $L'_n$ the lengths
of the curves $T^{-n}_1{\cal F}_{x,y}$ and $T^{-n}_1{\cal
F}'_{x,y'}$ respectively, where ${\cal F}_{x,y}$ is the part of
the curve in ${\cal F}_{x}$ between $x$ and $y$, and ${\cal
F}'_{x,y'}$ is understood in a similar way. By the assumption
stated in Part (c), the ratio between $A_n$ and $L_n\cdot L_n'$
tends to infinity, a contradiction.
\qed

\begin{Lemma}\label{lemma4.1}
For any $v, v'\in {\cal C}_x$,
$$
\sin\angle(DT_x(v), DT_x(v'))\le (1-d|x|^\c)\sin\angle(v, v').
$$
where the symbol $\angle(v,v')$ denotes the angle between the
vectors $v$ and $v'$.
\end{Lemma}

\proof Note that
$$|\det DT_x|_{E(v, v')}|
= {|DT_x(v)|\cdot |DT_x(v')|\cdot \sin\angle(DT_x(v),  DT_x(v'))
   \over |v|\cdot |v'|\cdot \sin\angle(v, v)}
$$
and
$$||DT_x|_{E(v)}|| = {|DT_x(v)|\over |v|}, \qquad
  || DT_x|_{E(v')}||= {|DT_x(v)|\over |v'|}.
$$
Then the results follows from (\ref{B.strong}).
\qed

\begin{Lemma}\label{lemma4.2}
There exist constants $H>0$, $a>0$, and
$\tau_1\in (0,1)$, such that for all $x\in TR\backslash \{p\}$,
\begin{equation}\label{L4,2}
d(E_x, E_y)\le {H d(x,y)^{\tau_1}\over |x|^{\tau_1} } \qquad
  \forall y\in B(x, a|x|),
\end{equation}
where $d(E_x, E_y)$ is defined by $d(E_x, E_y)=\sin\angle(v_x, v_y)$,
$v_x$ and $v_y$ are the tangent vectors of ${\cal F}_x$
and  ${\cal F}_y$ at $x$ and $y$ respectively chosen
in the way that $0\le \angle(v_x, v_y)<\pi/2$.
\end{Lemma}

\proof We note that we only need prove (\ref{L4,2}) for all $x$
in a small neighborhood $\tilde R\subset R$ of $p$, because
$DT_x(E_x)=E_{Tx}$, and then the results can be extended to $TR$.

Take $\tilde d\in (0,d)$.  Then for each $x$ we can extend ${\cal
C}_x$ to $\tilde {\cal C}_x$ such that (\ref{B.strong}) hold with
$\tilde d$ for all $v\in {\cal C}_x$ and $v'\in \tilde {\cal
C}_x$.  By (\ref{B.0}) and the fact that $T$ is $C^{1+\c}$, we
can write $DT(x)=T_0(x)+T_\c(x)+T_h(x)$, where $T_0=DT_p$, $T_\c$
satisfies $T_\c(tx)=t^\c T_\c(x)$\ $\forall t>0$ and
$|T_h(x)|=o(|x|^{\c})$.  So it is easy to see that we
can find $\e_a>0$ such that $\tilde {\cal C}_x\cap {\Bbb
S}^{m-1}$ contains an $\e_a$-neighborhood of ${\cal C}_x\cap
{\Bbb S}^{m-1}$ in ${\Bbb S}^{m-1}$ for all $x$ with $|x|$ small.
Moreover, since ${\cal C}_x$ is uniformly continuous in $(t,
\phi)$, we can take $a>0$ and $\tilde R$ small such that for all
$x\in \tilde R$, with $d(x,y)\le a|x|^\c$, ${\cal C}_y\subset
\tilde {\cal C}_x$. So if $v\in {\cal C}_x$ and $v'\in {\cal
C}_y$, we have
\begin{eqnarray*}\label{B.stronga}
{|\det DT_x|_{E(v, v')}|\over
\|DT_x|_{E(v)}\|\cdot \|DT_x|_{E(v')}\|}
\le 1-\tilde d|x|^\c.
\end{eqnarray*}
Hence, by the same arguments used in Lemma~\ref{lemma4.1} we have
\begin{eqnarray}\label{angle}
\sin\angle(DT_x(v), DT_x(v'))\le (1-\tilde d|x|^\c)\sin\angle(v, v').
\end{eqnarray}

Take $\tau_1\in (0,1)$ such that
\begin{equation}\label{tau1}
\Bigl(1-{\tilde d\over 2}|x|^\c\Bigr)
\Bigl({|Tx|\over |x|}\cdot {d(x, y)\over d(Tx,Ty)}\Bigr)^{\tau_1}
\le 1
\end{equation}
for all $x\in \tilde R$ close to $p$ with $d(x,y)\le a|x|$.

Take $0<a_1\le a$ such that if $d(x,y)\le a_1|x|$, then
\begin{equation}\label{barC2}
\|DT(x)-DT(y)\|\le \bar C_2|x|^{\c-1}d(x,y)^{\tau_1}.
\end{equation}
for some $\bar C_2>0$.  This is possible because of (\ref{B.2}).

Take $H>0$ such that $H\tilde d>2\bar C_2$.

Let ${\cal L}=\{L_x: x\in \tilde R\setminus \{p\}\ \}$ be the set
of all line bundles in the tangent bundle over $\tilde R$.
Clearly $DT$ induces a map ${\cal D}: {\cal L}\to {\cal L}$
given by $({\cal D}{\cal L})_x=DT_x(L_{T^{-1}_1x})$,
and ${\cal E}=\{E_x\}$ is the unique fixed point of ${\cal D}$
contained in ${\cal C}$.  Denote
\begin{equation}\label{spaceH}
{\cal H}=\Bigl\{ \{L_x\}\in {\cal L}\cap {\cal C}:
d(L_x, L_y)\le {H d(x,y)^{\tau_1}\over |x|^{\tau_1} } \qquad
  \forall y\in B(x, a_1|x|)\ \Bigr\}.
\end{equation}
We show that ${\cal D}({\cal H})\subset {\cal H}$.
This implies the result since
$\{E_x\}=\cap_{n\ge 0}{\cal D}^n{\cal C}$.

Take $\{L_x\}\in {\cal H}$. Let $x,y\in \tilde R$ with $d(x,y)\le
a_1|x|$. Take unit vectors $e_x\in L_x$, $e_y\in L_y$. So
$\sin\angle(e_x,e_y)\le H|x|^{-\tau_1}d(x,y)^{\tau_1}$. By
(\ref{angle}) and (\ref{barC2}),
\begin{eqnarray*}
  &&  \sin\angle(DT_x(e_x), DT_y(e_y))  \\
&\le &\sin\angle(DT_x(e_x), DT_x(e_y))
   + \sin\angle(DT_x(e_y), DT_y(e_y))  \\
&\le &(1-\tilde d|x|^\c) \sin\angle(e_x, e_y)+ |DT_x(e_y)-DT_y(e_y)|  \\
&\le &(1-\tilde d|x|^\c){H d(x,y)^{\tau_1}\over |x|^{\tau_1} }
   +  {\bar C}_2|x|^{\c-1}d(x,y)^{\tau_1}   \\
&=& \bigl[(1-\tilde d|x|^\c)H+{\bar C}_2|x|^\c\bigr]
  {d(Tx, Ty)^{\tau_1}\over |Tx|^{\tau_1}}
 \cdot {d(x,y)^{\tau_1}
  \over d(Tx, Ty)^{\tau_1}} {|Tx|^{\tau_1}\over |x|^{\tau_1}}.
\end{eqnarray*}
By the choice of $H$, the quantity in the blanket is less
than $1-\tilde d|x|^\c/2$.  Then by (\ref{tau1}) the right side of the
inequality is less than $H|Tx|^{-\tau_1}d(Tx,Ty)^{\tau_1}$.
We get the desired results.
\qed

\begin{Lemma}\label{lemma4.3}
There exists $J^*>0$ such that for any $x,y$ with
$d(x_i,y_i)\le |x_i|^{\bar \c}$ for some $\bar \c>1$,
$i=1,\cdots, n$,
\begin{equation}\label{distF}
{|DT_1^{-n}(y)|_{E_y({\cal F})}|
   \over |DT_1^{-n}(x)|_{E_x({\cal F})}|}
\le J^*.
\end{equation}
\end{Lemma}

\proof
Take an integer $\bar r\ge 2C_0'/C_0$,
where $C_0$ and $C_0'$ are as in (\ref{B.1}).
We assume that
$x_0\le 1/(\c  C_0' k_0)^\b$
for some $k_0\ge 1$.
Then we take $k_i =(\bar r^i -1)k_0$ for $i =1, \cdots, \ell-1$,
where $\ell-1$ is the largest number $j$ such that $k_j< n$.
Let $k_\ell=n$.
By Lemma~\ref{lemma3.1}, we know that
\begin{eqnarray}\label{xj}
|x_{j}|^\c \le 2/(\c  C_0'  (k_0+j)).
\end{eqnarray}
Hence, (\ref{B.1}) implies
$$
\|DT^{k_i-k_{i-1}}_{x_{k_{i}}}\|
\le\!\!\prod_{j=k_{i-1}}^{k_i-1} \|DT_{x_j}\|
\le\!\! \prod_{j=k_{i-1}}^{k_i-1} \Bigl(1+{2C_1\over\c  C_0'(k_0+j)}\Bigr)
\le\!\! \prod_{j=k_{i-1}}^{k_i-1} \Bigl(1+{1\over k_0+j}\Bigr)^C
$$
for some $C$ larger than $2C/\c C_0'$ if $k_i$ is large enough.
So the choice of $\bar r$ gives
\begin{eqnarray}\label{ineqDT}
\|DT^{k_i-k_{i-1}}_{x_{k_{i}}}\|
\le \Bigl({k_0+k_i\over k_0+k_{i-1}}\Bigr)^C
\le \bar r^C
\end{eqnarray}
for all $i\ge 0$.


Let $e_x$ be the unit tangent vector of ${\cal F}$ at $x$.
We have
\begin{eqnarray*}\label{ratio}  \disp
  {|DT_1^{-n}(y)|_{E_y({\cal F})}|
   \over |DT_1^{-n}(x)|_{E_x({\cal F})}|}
&=& {|DT^{n}_{x_n}(e_{x_n})|  \over |DT^{n}_{y_n}(e_{y_n})|}
={|DT^{n}_{x_n}(e_{x_n})|  \over |DT^{n}_{x_n}(e_{y_n})|}
\cdot {|DT^{n}_{x_n}(e_{y_n})|  \over |DT^{n}_{y_n}(e_{y_n})|} \\
&=&\prod_{i=1}^{\ell}
  {|DT^{k_i-k_{i-1}}_{x_{k_{i}}}(e_{x_{k_{i}}})|
      \over |DT^{k_i-k_{i-1}}_{x_{k_{i}}}(e_{y_{k_{i}}})|}
  \cdot \prod_{j=1}^{n}  {|DT_{x_j}(e_{y_j})|
           \over |DT_{y_j}(e_{y_j})|}.
\end{eqnarray*}

By the results of Lemma~\ref{lemma4.2} and (\ref{ineqDT}),
each factor in the first product is bounded by
\begin{eqnarray*}\label{sin}
&&1+{|DT^{k_i-k_{i-1}}_{x_{k_{i}}}(e_{x_{k_{i}}})|
                -|DT^{k_i-k_{i-1}}_{x_{k_{i}}}(e_{y_{k_{i}}})|
           \over |DT^{k_i-k_{i-1}}_{x_{k_{i}}}(e_{y_{k_{i}}})|}  \\
&\le& 1+{|DT^{k_i-k_{i-1}}_{x_{k_{i}}}(e_{x_{k_{i}}}-e_{y_{k_{i}}})|
           \over |DT^{k_i-k_{i-1}}_{x_{k_{i}}}(e_{y_{k_{i}}})|}
\le 1+{\|DT^{k_i-k_{i-1}}_{x_{k_{i}}}\|\cdot |e_{x_{k_{i}}}-e_{y_{k_{i}}}|
          \over |DT^{k_i-k_{i-1}}_{x_{k_{i}}}(e_{y_{k_{i}}})|}   \\
&\le& 1+{\bar r^C \cdot BHd(x_{k_{i}}, y_{k_{i}})^{\tau_1}
          \over |x_{k_i}|}
\le 1+\bar r^C BH|x_{k_i}|^{\tau_1(\bar \c-1)},
\end{eqnarray*}
where we use the fact that $|e_{x_{k_{i}}}-e_{y_{k_{i}}}|
\le B\sin\angle(e_{x_{k_{i}}}, e_{y_{k_{i}}})$ for some $B>0$.
Also note that by (\ref{xj}) and the choice of $k_i$,
$\{|x_{k_i}|\}$ decreases exponentially fast as $i\to\infty$.
Since $\bar \c>1$, the first product in above equality is convergent.

For the second product, by (\ref{B.2}) each factor is bounded by
\begin{eqnarray*}
&&1+{|DT_{x_j}(e_{y_j})|-|DT_{y_j}(e_{y_j})|
           \over |DT_{y_j}(e_{y_j})|}
\le 1+{C|x_j|^{\c-1}d(x,y)  \over |DT_{y_j}(e_{y_j})|}
\le 1+{C|x_j|^{\bar\c+\c-1} \over |DT_{y_j}(e_{y_j})|}.
\end{eqnarray*}
By (\ref{xj}) and the fact $\bar\c>1$,
we know that $\sum_j |x_j|^{\bar\c+\c-1}$ converges.
So the product is also bounded.  We get the result.
\qed

\section{Proof of Theorem A}
\setcounter{equation}{0}

In this section we first introduce a subspace $V_\a$ of $L^1\equiv
L^1(\Bbb R^m, \nu)$  with compact unit ball that contains the
density function of the invariant measures of the induced map of $T$
with respect to the relatively compact subspace $M\backslash R$. Here
we only give a brief description and list some properties we use.
We refer to \cite{S} and \cite{K} for more details.

Let $f$ be an $L^1(\Bbb R^m, \nu)$ function.  If $\Omega$ is a Borel
subset of $\Bbb R^m$, we define the oscillation of $f$ over
$\Omega$ by the difference of essential supremum and essential
infimum of $f$ over $\Omega$:
$$
\osc(f,\Omega) = \Esup_\Omega f - \Einf_\Omega f.
$$
If $B_{\epsilon}(x)$ denotes the ball of radius $\epsilon$ about
the point $x$, then we get a measurable function $x\rightarrow
\mbox{osc}(f, \ B_{\epsilon}(x))$.  The function have the
following properties.

\begin{Proposition}\label{properties1}
  Let $f, f_i, g\in L^\infty(\Bbb R^m,\nu)$ with $g\ge 0$, $\e>0$, and $S$
  be a Borel subset of $\Bbb R^m$.  Then
\begin{enumerate}
\item[{\rm (i)}]
$\disp \osc\bigl(\sum_i f_i,  \ B_\e(\cdot)\bigr)
\le \sum_i \osc\bigl(f_i, \ B_\e(\cdot)\bigr)$,

\item[{\rm (ii)}]
$\disp \osc\bigl(f {\Bbba}_S,  \ B_\e(\cdot)\bigr)
 \le \osc\bigl(f, \ S\cap B_\e(\cdot)\bigr){\Bbba}_S(\cdot)
+ 2\bigl[\Esup_{S\cap B_\e(\cdot)} f \bigr]\Bbba_{B_\e(S)\cap B_\e(S^c)}$,

\item[{\rm (iii)}]
$\disp \osc\bigl(fg, \ S\bigl)
\le \osc\bigl(f, \ S\bigl)\Esup_S g
   +\osc\bigl(g, \ S\bigl)\Einf_S f$.
\end{enumerate}
\end{Proposition}

\proof See \cite{S} Proposition 3.2. \qed

Take $0<\a <1$ and $\e_0>0$.
We define the $\alpha$-seminorm of $f$ as:
\begin{eqnarray}\label{seminorm}
|f|_{\alpha} = \sup_{0<\epsilon\le\epsilon_0}
\epsilon^{-\alpha} \int_{\Bbb R^m} \mbox{osc}(f,
B_{\epsilon}(x)) d\nu(x).
\end{eqnarray}
We will consider the space of the functions $f$
with bounded $\alpha$-seminorm, namely:
\begin{eqnarray}\label{norm}
V_{\alpha} =  \left\{ f\in L^1: |f|_{\alpha} <\infty \right\}
\end{eqnarray}
and equip $V_{\alpha}$ with the norm:
\begin{eqnarray}\label{falphanorm}
\parallel\cdot\parallel_{\alpha} =
\parallel\cdot\parallel_1 + |\cdot|_{\alpha},
\end{eqnarray}
where $\|\cdot\|_1$ denotes the $L^1$ norm. This space will not
depend on the choice of $\epsilon_0$.  With the
$\parallel\cdot\parallel_{\alpha}$ norm, $V_\a$ is a Banach
space; moreover according to Theorem~1.13 in \cite{K}, the unit ball
in $V_\a$ is compact in $L^1$.

\begin{Proposition}\label{properties2}
Let $f\in V_\a$; then:
\begin{enumerate}
\item[{\rm (i)}]
$\disp \|f\|_\infty \le {1\over \c_m\e_0^m}\|f\|_\a$
provided $\e_0\le 1$.

\item[{\rm (ii)}]
There exists a ball $B_\e(x)$ such that
$\disp \Einf_{B_\e(x)} f>0$.
\end{enumerate}
\end{Proposition}

\proof See \cite{S} Proposition 3.4 and Lemma 3.1. \qed

To prove Theorem A  we need one more ingredient, the so-called
Lasota-Yorke's inequality, which will be proved in Section 6. This
inequality provides an upper bound on the action of the
Perron-Frobenius operator on the elements on $V_\a$. Such an
operator will be defined on the subspace $M\backslash R$ with a
potential given by the inverse of the determinant of the induced
map. We will denote it as $\hat Pf$. We will prove that:
$$
|\hat Pf|_\a\le \eta|f|_\a+D\|f\|_1
$$
where $\eta<1$ and $D<\infty$. This, plus the compactness in $L^1$
of the unit ball of $V_\a$, will allow us to invoke the ergodic
theorem of Ionescu-Tulcea and Marinescu \cite{ITM} (see also
\cite{K}, Theorem 3.3), to conclude that there exists a unique
(greatest)\footnote {``Unique greatest'' means that any other
measure absolutely continuous with respect to $\nu$ is absolutely
continuous with respect to $\mu$.} invariant probability measures
$\mu$ which is absolutely continuous with respect to $\nu$ on
$M\backslash R$ and which decomposes into a finite number of
cyclic disjoint measurable sets upon which a certain power of the
map is mixing.

\medskip
\proofa

Recall that $R$ is given in Assumption 3. We construct a induced
system $(\hat M, \hat T)$. Denote $\hat M=M\backslash R$. Let
$\hat T: \hat M\to \hat M$ be the first return map of $T$, so
that  $\hat T(x)=T(x)$ if $x\not\in  T^{-1}R$, otherwise
$\hat T(x)=T^{i+1}(x)=T_1^iT_j(x)$ if $x\in T_j^{-1}R$, where $i$ is
the smallest positive integer such that $T_1^iT_j(x)\notin R$.
We denote $g(x)=|\det DT(x)|^{-1}$, and
similarly $\hat g(x)=g(x)$ if $x\not\in  T^{-1}R$
and $\hat g(x)=|\det DT^{i+1}(x)|^{-1}$ if otherwise.
Let $\hat \m$
be the conditional measure of the Lebesgue measure $\m$. We may
still think that $\hat \m$ is a Lebesgue measure with $\hat
\m(\hat M)=1$.

Let $P$ be the Perron-Frobenius operator of $T$ with the potential
function $\log g(x)$, i.e.
$$
Pf(x)=\sum_{Ty=x} f(y)g(y).
$$
Then let $\hat P$ be the Perron-Frobenius operator of $\hat T$
with the potential function $\log \hat g(x)$.

By Proposition~\ref{LYineq} in the next section we have the
Lasota-Yorke's inequality for the induced system $(\hat M, \hat
T)$. So $\hat T$ has an absolutely continuous invariant
probability measure $\hat \mu$ on $\hat M$ with density function
$\hat h$ that has finitely many ergodic components.

We extend $\hat \mu$ to $M$ to get an invariant measure of $T$.
Recall $R_0=TR\backslash R$,  and let $R_n=T_1^{-n}R_0$ for $n>0$.
By Remark~\ref{Rmk3a}, $\diam R_n\to 0$.
So we have $R=\sum_{n=1}^\infty R_n\cup \{p\}$.
We extend $\hat h$ to $R$ to get a density function $h$ on $M$.
That is, if $h$ is defined on $M\backslash T_1^{-n}R$, then for
$x\in T_1^{-n}R\backslash T_1^{-n-1}R$, we let
$$
h(x)= g(x)^{-1}\cdot \Bigl(h(Tx)
     -\sum_{j\not= 1} h(T_j^{-1}Tx)g(T_j^{-1}Tx)\Bigr).
$$
It is easy to see that $h\ge 0$ and $Ph=h$ on $M$.
Let $\mu$ be the measure on $M$ with density $h$.
Clearly, $\mu$ is invariant under $T$ and has the same number of
ergodic components as $\hat \mu$ does.

Next, we show that $\mu M$ is finite if
$\disp \sum_{i=1}^\infty \nu(T_1^{-i} R)<\infty$.
Since $\mu$ is invariant, we have
$$
\mu R_i=\mu R_{i+1} + \sum_{j=2}^{K'} \mu(T_j^{-1}R_i),
$$
where we assume that in addition to $T_1^{-1}R\subset U_1$,
$R$ has $K'-1$ preimages in $U_2, \cdots, U_{K'}$, where $K'\le K$.
Take summation from $i=n$ to infinity, we get
$$
\mu R_n=\sum_{j=2}^{K'} \mu\bigl(T_j^{-1}\bigcup_{i=n}^\infty R_i\bigr)
=\sum_{j=2}^{K'} \mu\bigl(T_j^{-1}T_1^{-n}R\bigr).
$$
Note that $\|\hat h\|_\infty\le \infty$ since $\hat h\in V_\a$, and
then note that the Jacobian of $T_j^{-1}$ is less than or equal to $1$.
We have
$$
    \mu\bigl(T_j^{-1}T_1^{-n}R\bigr)
\le \|\hat h\|_\infty \m\bigl(T_j^{-1}T_1^{-n}R\bigr)
\le \|\hat h\|_\infty \m\bigl(T_1^{-n}R\bigr).
$$
Hence
\begin{eqnarray}\label{upperbound}
\mu R =\sum_{n=1}^\infty \mu R_n
\le \|\hat h\|_\infty (K'-1) \sum_{n=1}^\infty \m\bigl(T_1^{-n}R\bigr)
 < \infty.
\end{eqnarray}

Now we prove the last part of the theorem. By
Proposition~\ref{properties2}(ii), there is a ball
$B_\e(z)\subset M\setminus R$ such that $\disp \Einf_{B_\e(x)}
\hat h \ge h_*>0$ for some constant $h_*$. By our assumption,
there exists $\tilde N>0$ such that $T^{\tilde N}B_\e(z)\supset
M$. So for any $x\in M$, there is $y_0\in B_\e(z)$ such that
$T^{\tilde N}y_0=x$. Since $|\det DT|$ is bounded above, we have
$g_*:=\inf\{g(y): \ y\in M\}>0$. Hence, for  every $x$,
$$
h(x)=(P^{\tilde N}h)(x)
=\sum_{T^{\tilde N}y=x} h(y)\prod_{i=0}^{{\tilde N}-1} g(T^iy)
\ge h(y_0)\prod_{i=0}^{{\tilde N}-1} g(T^iy_0)
\ge h_* g_*^{\tilde N}.
$$
In this case, we can use a similar method as for (\ref{upperbound})
to get
\begin{eqnarray*}
\mu R =\sum_{n=1}^\infty \mu R_n
\ge (h_* g_*^{\tilde N}) g_* \sum_{j=2}^{K'} \m\bigl(T_1^{-n}R\bigr)
= \infty.
\end{eqnarray*}
This ends the proof.
\qed

\section{A Lasota-Yorke type inequality}
\setcounter{equation}{0}

Let $R$ be as in Assumption 3.
Denote $\hat T_{ij}=T_1^iT_j$ and
$U_{ij}=\hat T_{ij}^{-1} (R_0)= T_j^{-1}R_i$ for $i>1$ and
$U_{0j}=U_j\backslash T_j^{-1} R$.  So if $TU_l\not\ni p$,
then $U_{il}$ is undefined for any $i>0$ and $U_{0l}=U_l$.
Clearly, $U_{ij}\subset U_j$ for all $i>0$ and $\{U_{ij}, i\ge0\}$ are
pairwise disjoint.

\begin{Lemma}\label{proposition1}
There exists $0<\e_6\le \e_5$ such that
for any $\e_0\le \e_6$, $\e\le \e_0$, $x\in M$,
\begin{eqnarray}\label{prop1}
2\sum_{j=1}^K \sum_{i=0}^{\infty} {\m(\hat
T_{ij}^{-1}B_\e(\partial R_0) \cap B_{(1-s)\e_0}(x)) \over
\m(B_{(1-s)\e_0}(x))} \le {\lambda \e^\a \over \e_0^\a},
\end{eqnarray}
where $\l$ is given by Assumption 3(b).
\end{Lemma}

\proof
Note that the sets
$\cup_{i=1}^\infty \partial U_{ij}$, $j=1,\cdots, K$,
are pairwise separated.  So by Assumption 3(b) and the definition
$\lambda$ in (\ref{lambda}) we only need
prove that there exists $\e_6>0$ such that for any given $j$,
for any $x$ in the $\e_6$-neighbourhood of $T_j^{-1}R_0$,
if $0<\e\le \e_0\le \e_6$, then
\begin{eqnarray}\label{2sum}
2\sum_{i=0}^{\infty} {\m(\hat T_{ij}^{-1}B_\e(\partial R_0) \cap
B_{(1-s)\e_6}(x)) \over \m(B_{(1-s)\e_6}(x))} \le {\lambda \e^\a
\over \e_0^\a}.
\end{eqnarray}

Take
\begin{eqnarray*}\label{ep5}
\e_6\le \min\{\e_5, \e_3\}\cdot
\Bigl({\l (1-s)^m \over 2C_\xi I^2}\Bigr)^{1/\a},
\end{eqnarray*}
where $\e_3$ is given by Assumption 3(c).

Recall that $N_s$ is also given by Assumption 3(c).
Reduce $\e_6$ if necessary
such that for any $x$, the ball $B_{(1-s)\e_6}(x)$
intersects at most one connected component of the set
$\{\hat T_{ij}^{-1} B_{\e_6}(\partial R_0), 0<i\le N, 1< j\le K\}$.
We also require $\e_6$  small enough such that
for any $1< j\le K$, $1\le i\le N_s$,
the part $\hat T_{ij}^{-1} \partial R_0\cap B_{\e_6}(x)$
are close to an $(m-1)$ dimensional plane.

Take $\e$ and $\e_0$ such that $0<\e\le \e_0\le \e_6$ .

We first consider the case $1\le i\le N_s$.
Note that
$\hat T_{ij}^{-1}  B_{\e}(\partial R_0)\cap B_{(1-s)\e_0}(x)
\subset
B_{s\e}(T_{ij}^{-1} \partial R_0)\cap B_{(1-s)\e_0}(x)$.
The volume of the latter is close to
$\gamma_{m-1}((1-s)\e_0)^{m-1}\cdot 2s\e
=2s\gamma_{m-1}\e(1-s)^{m-1}\e_0^{m-1}$.
So
$\disp {\m(\hat T_{ij}^{-1}B_\e(\partial R)\cap B_{(1-s)\e_0}(x))
\over \m(B_{(1-s)\e_0}(x))}$
is close to
$\disp {2s\gamma_{m-1}\e(1-s)^{m-1}\e_0^{m-1}
\over \gamma_{m}(1-s)^{m}\e_0^{m}}
={2s\gamma_{m-1}\e\over (1-s)\gamma_{m}\e_0}$.
Hence, by Assumption 3(b), we know that it is less than
$\lambda\e^\a/\e_0^\a$.

Now we consider the case that $i\ge N_s$.


Let
$\disp {\tilde \e}=\e_0\Bigl({2C_\xi I^2\over \l (1-s)^m}\Bigr)^{1/\a}$.
we have ${\tilde \e}\le \e_5$.

For each $i$, we take a partition
$\xi_i=\{\tilde A_{i1}, \tilde A_{i2}, \cdots, \}$
satisfying Assumption 4(c) with $n=i$ and $\tilde\e\le \e_5$.
Denote
$A_{ik}=\tilde A_{ik}\cap B_{{\tilde \e}}(\partial R_0)$,
$A_{ik}'=\tilde A_{ik}\cap B_{{\e}}(\partial R_0)$,
$A_{ijk}=\hat T_{ij}^{-1}A_{ik}$ and
$A_{ijk}'=\hat T_{ij}^{-1}A_{ik}'$.
Then we let
$$
{\cal A}=\{A_{ijk}: \ A_{ijk}'\cap B_{(1-s)\e_0}(x)\not=\emptyset\},
\qquad
{\cal A'}=\{A_{ijk}': \ A_{ijk}\in \cal A\}.
$$
By abusing notations, we may also think that ${\cal A}$
and ${\cal A}'$ are the unions of the sets they contain.

By the fact
$$
\nu A_{ijk}=\int_{A_{ik}} |\det D{\hat T}^{-1}_{ij}(x)| d\nu(x)
$$
and Assumption 4(c), we know that
\begin{eqnarray}\label{nu1}
  {\nu A_{ijk}'\over \nu A_{ijk}}
\le {C_\xi \e^\a \over {\tilde \e}^\a}\cdot I^2
=   {C_\xi I^2\e^\a \l (1-s)^m\over 2C_\xi I^2\e_0^\a}
=   {\e^\a \l (1-s)^m\over 2\e_0^\a}.
\end{eqnarray}

Denote
$s^*=\sup\bigl\{s(T_1^{-N_s}(z), T_1^{N_s}):
    \ z\in B_{\tilde \e}(R_0)\bigr\}$.
Note that by Assumption 4(c),
$\disp \diam A_{ik}\le 5m{\tilde \e}
\le 5m \e_0\Bigl({2C_\xi I^2\over \l (1-s)^m}\Bigr)^{1/\a}$.
Since $i\ge N_s$,  by Assumption 3(c), we have
\begin{eqnarray}\label{diamAijk}
\diam A_{ijk}
\le 5m \e_0\Bigl({2C_\xi I^2\over \l (1-s)^m}\Bigr)^{1/\a}\cdot s^*
=  s\e_0.
\end{eqnarray}
So if $A_{ijk}\in {\cal A}$, then
$A_{ijk}\cap B_{(1-s)\e_0}(x)\not=\emptyset$,
and therefore
$A_{ijk}\subset B_{\e_0}(x)$.  That is,
\begin{eqnarray}\label{nu2}
{\cal A}\subset B_{\e_0}(x).
\end{eqnarray}
Note that
\begin{eqnarray}\label{nu3}
\bigcup_{i=0}^{\infty}
\hat T_{ij}^{-1}B_\e(\partial R)\cap B_{(1-s)\e_0}(x)
\subset {\cal A}'.
\end{eqnarray}
By (\ref{nu1})-(\ref{nu3}), we get
\begin{eqnarray*}\label{result3.0}
&&2\sum_{i=0}^{\infty} {\m(\hat T_{ij}^{-1}B_\e(\partial R)
\cap B_{(1-s)\e_0}(x)) \over \m(B_{(1-s)\e_0}(x))}
\le 2\cdot  {\nu {\cal A}'\over \nu {\cal A}}
 \cdot {\nu {\cal A} \over \mu B_{{\e_0}}(x)}
 \cdot {\mu B_{{\e_0}}(x) \over \m(B_{(1-s)\e_0}(x)} \\
&\le& 2\cdot   {\e^\a \l (1-s)^m\over 2\e_0^\a}  \cdot 1
 \cdot {\c_m\e_0^m \over \c_m(1-s)^m\e_0^m}
 =\l {\e^\a\over \e_0^\a}.
\end{eqnarray*}
This is (\ref{2sum}), the formula we need show.
\qed


\begin{Proposition}\label{LYineq}
Assume that $T: M\to M$ satisfies Assumption 1-4,
and $\hat T: \hat M\to \hat M$ is the reduced system with respect
to $\hat M=M\backslash R$.
Then there exist $\eta < 1$ and $D< \infty$ such that for any
$f\in V_\a=V_\a(\e_0)$, we have $Pf\in V_\a$ and
$$
|\hat Pf|_\a\le \eta|f|_\a+D\|f\|_1
$$
for all $\e_0$ sufficiently small.
\end{Proposition}

\proof
Take $\dd>0$ such that for any $\e\le \e_4$,
\begin{eqnarray}\label{chiep}
(1+Js^\a \e^a)(1+cs^a\e^\a)\le 1+\dd\e^\a,
\end{eqnarray}
where $c$ and $J$, $\e_4$ and are given in Assumption 4(a) and (b)
respectively.

Recall that by Assumption 3(b), $s^\a+\lambda\le \eta_0 <1$.
Take $b>0$ such that $(s^\a+ \lambda)+3K'b\c_m^{-1}< 1$,
where $K'$ is the number of preimages of $p$ for the map $T$.
Recall also that $\e_1, \e_2$, $\e_4$ and $\e_6$ are given in Assumption 1(b),
3(b) and 4(b) and Lemma~\ref{proposition1} respectively.
Take $\e_0\le \min\{\e_1, \e_2, \e_4, \e_6\}$ such that
\begin{eqnarray}\label{xib}
\eta := (1+\dd \e_0^\a)(s^\a+ \lambda)+3K'b\c_m^{-1}< 1.
\end{eqnarray}


Denote
$$
G_R(x, \e,\e_0) =2\sum_{j=1}^K \sum_{i=0}^{N(\e)} {\m(\hat
T_{ij}^{-1}B_\e(\partial R_0) \cap B_{(1-s)\e_0}(x)) \over
\m(B_{(1-s)\e_0}(x))}.
$$
Recall that $G_U(x, \e,\e_0)$ is given by (\ref{GU})
in Assumption 3(b).
Note that if $\e_0$ is small, then $\supp G_U(\cdot, \e,\e_0)$ and
$\supp G_R(\cdot, \e,\e_0)$ are disjoint.
 Also, by (\ref{lambda})
and Lemma~\ref{proposition1}, we know that
\begin{eqnarray}\label{G<lambda}
G(\e,\e_0)=\sup_{x\in M}\{G_U(x, \e,\e_0), G_R(x, \e,\e_0)\}
\le {\l\e^\a\over \e_0^\a}.
\end{eqnarray}
Then we take
\begin{eqnarray}\label{DD}
D = 2\dd+ 2(1+\dd \e^\a)\sup_{\e\le \e_0}G(\e, \e_0)\e^{-\a}
   +K'b\c_m^{-1}.
\end{eqnarray}
By (\ref{G<lambda}), $G(\e, \e_0)\e^{-\a}\le \l \e_0^{-\a}$.
We have $D<\infty$.

Let $\e\le \e_0$.

By Proposition~\ref{properties1},
\begin{eqnarray}\label{osc(pf)}
 &&\osc\bigl(\hat Pf, \ B_\e(x)\bigr)
\le \sum_{j=1}^K \sum_{i=0}^\infty
  \osc\bigl((f {\hat g})\circ \hat T_{ij}^{-1}\Bbba_{\hat TU_{ij}}, \;
           B_\e(x)\bigr) \nonumber\\
&\!\!\!\!\!\!\le&\!\!\!\!\!\! \sum_{j=1}^K \sum_{i=0}^\infty \Bigl(
  \osc\bigl((f {\hat g})\circ \hat T_{ij}^{-1},  B_\e(x)\bigr)
        \Bbba_{\hat TU_{ij}}(x)
  +\bigl[2\Esup_{B_\e(x)}(f {\hat g})\circ \hat T_{ij}^{-1} \bigr]
         \Bbba_{B_\e(\partial \hat TU_{ij})}(x)
   \Bigr)                                                 \nonumber\\
&\!\!\!\!\!\!=:&\!\!\!\!\!\! \sum_{j=1}^K \sum_{i=0}^\infty \Bigl(
  R^{(1)}_{ij}(x)\Bbba_{\hat TU_{ij}}(x)
 +R^{(2)}_{ij}(x)\Bbba_{B_\e(\partial \hat TU_{ij})}(x)  \Bigr).
\end{eqnarray}

Denote $y_{ij}=\hat T_{ij}^{-1}x$.
We can choose $N=N(\e)>0$ for each $0<\e\le \e_0$
according to Assumption 4(b).

For $R^{(1)}_{ij}(x)$ with $x\in \hat TU_{ij}$,
we first consider the case $i\le N(\e)$.
By Assumption 4(a), (b) and (\ref{chiep}), we have
${\hat g}(y_{ij}')/{\hat g}(y_{ij})
\le (1+Js^\a\e^\a)(1+cs^\a\e^\a)\le 1+\dd \e^\a$
if $d(T^{i+1}y_{ij}, T^{i+1}y_{ij}')\le s\e$.  Hence
${\hat g}(y_{ij}')\le (1+\dd \e^\a){\hat g}(y_{ij})$ and
$\osc\bigl({\hat g}, \ B_{s\e}(y_{ij}))\le 2\dd \e^\a{\hat g}(y_{ij})$.
So we get
\begin{eqnarray*}\label{est2}
 && R^{(1)}_{ij}(x)
= \osc\bigl(f {\hat g}, \ \hat T_{ij}^{-1} B_\e(x)\cap U_{ij}\bigr)
                                                           \nonumber\\
&\le&\!\!\!\osc\bigl(f, \ B_{s\e}(y_{ij}) \cap U_{ij}\bigr)
   \Esup_{B_{s\e}(y_{ij})\cap U_{ij}} {\hat g}
 + \osc\bigl({\hat g}, \ B_{s\e}(y_{ij})\cap U_{ij}\bigr)
   \Einf_{B_{s\e}(y_{ij})\cap U_{ij}} f                    \nonumber\\\
&\le& \!\!\!(1+\dd \e^\a)\osc\bigl(f, \ B_{s\e}(y_{ij}) \cap U_{ij}\bigr)
     {\hat g}(y_{ij})
  +2\dd \e^\a f(y_{ij}){\hat g}(y_{ij}).
\end{eqnarray*}
If $i>N(\e)$, then we must have $x\in R_0$, and therefore
\begin{eqnarray*}
 &&{{}}\!\!\!R^{(1)}_{ij}(x)
=\osc\bigl(f {\hat g}, \ \hat T_{ij}^{-1} B_\e(x)\cap U_{ij}\bigr)
                                                           \nonumber\\
&\le&\!\!\!\osc\bigl(f, \ B_{s\e}(y_{ij}) \cap U_{ij}\bigr)
   \Einf_{B_{s\e}(x)\cap U_{ij}} {\hat g}
 + \osc\bigl({\hat g}, \ B_{s\e}(y_{ij})\cap U_{ij}\bigr)
   \Esup_{\hat T_{ij}^{-1}B_{\e}(x)} f      \nonumber\\
&\le&\!\!\!\osc\bigl(f, \ \  B_{s\e}(y_{ij}) \cap U_{ij}\bigr) {\hat g}(y_{ij})
 + \|f\|_\infty \sup_{\hat T_{ij}^{-1}B_{\e}(x)} {\hat g}.
\end{eqnarray*}
By Assumption 4(b), for any $x\in R_0$,
$\disp \sum_{i=N}^\infty
     (\sup_{\hat T_{ij}^{-1}B_{\e}(x)} {\hat g})\le b\e^{m+\a}$.
Hence,
\begin{eqnarray*}
 &&\!\!\!\sum_{j=1}^K \sum_{i=0}^\infty R^{(1)}_{ij}(x)
    \Bbba_{\hat TU_{ij}}(x)
\le  K'b\e^{m+\a}   \|f\|_\infty \Bbba_{R_0}(x)  \nonumber\ \\
&+ &\!\!\!\sum_{j=1}^K \sum_{i=0}^\infty
    \Bigl( (1+\dd \e^\a)\osc\bigl(f, \ B_{s\e}(y_{ij}) \cap U_{ij}\bigr)
            {\hat g}(y_{ij})
  +2\dd \e^\a f(y_{ij}){\hat g}(y_{ij}) \Bigr)          \nonumber\\
&\le& \!\!\!  K'b\e^{m+\a}  \|f\|_\infty \Bbba_{R_0}(x)  \! +
(1+\dd \e^\a)\bigl[\hat P\osc\bigl(f, \ B_{s\e}(\cdot)\bigr)\bigr](x)\!
 + 2\dd \e^\a (\hat Pf)(x).
\end{eqnarray*}
Since $\int_{\hat M}  \hat Pf d\hat \m=\int_{\hat M}  f d\hat\m$
for any integrable function $f$,  we have
\begin{eqnarray}\label{R1}
&& \int_{\hat M} \sum_{j=1}^K \sum_{i=0}^\infty
     R^{(1)}_{ij} \Bbba_{\hat TU_{ij}}d\hat \m            \nonumber\\
&\le &  K'b\e^{m+\a}  \|f\|_\infty \hat \m R_0
+ (1+\dd \e^\a)\int_{\hat M}\osc\bigl(f, \ B_{s\e}(\cdot)\bigr)d\hat\m
+ 2\dd \e^\a \int_{\hat M}f d\hat \m            \nonumber\\
&\le &(1+\dd \e^\a)s^\a \e^\a |f|_\a +2\dd \e^\a\|f\|_1
 +K'b\e^{m+\a}  \|f\|_\infty \hat \m R_0.
\end{eqnarray}

As for $R^{(2)}_{ij}(x)$, if $i\le N(\e)$, then we have
$$\disp \Esup_{B_\e(x)}(f {\hat g})\circ \hat T_{ij}^{-1}
\le \bigl[\Esup_{B_{s\e}(y_{ij})}|f| \bigr] \hat g(y_{ij}) (1+\dd \e^\a).
$$
Hence by the same method as in \cite{S}, we get that
\begin{eqnarray*}
 \int_{\hat M}\sum_{j=1}^K \sum_{i=0}^{N(\e)}
  R^{(2)}_{ij} \Bbba_{B_\e(\partial \hat TU_{ij})} d\hat \m
\le  2(1+\dd \e^\a)G(\e, \e_0)\bigl(\e_0^\a |f|_\a + \|f\|_1).
\end{eqnarray*}
If $i\ge N(\e)$, then
$\disp \Esup_{B_\e(x)}(f {\hat g})\circ \hat T_{ij}^{-1}
\le \|f\|_\infty \sup_{\hat T_{ij}^{-1}B_{\e}(x)} {\hat g}$,
and
\begin{eqnarray*}
\sum_{j=1}^K \sum_{i=N(\e)}^\infty
  R^{(2)}_{ij} \Bbba_{B_\e(\partial \hat TU_{ij})}
\le 2K'\|f\|_\infty
  \sum_{i=N(\e)}^\infty \sup_{\hat T_{ij}^{-1}B_{\e}(x)} {\hat g}
\end{eqnarray*}
Again, by Assumption 4(b) it is bounded by $2K'b\e^{m+\a}\|f\|_\infty$.
So we have
\begin{eqnarray}\label{R2}
&& \int_{\hat M}\sum_{j=1}^K \sum_{i=0}^{\infty}
  R^{(2)}_{ij} \Bbba_{B_\e(\partial \hat TU_{ij})} d\hat \m
                           \nonumber\\
&\le&\!\!\!  2(1+\dd \e^\a)G(\e, \e_0)\bigl(\e_0^\a |f|_\a + \|f\|_1)
    +2K'b\e^{m+\a}\|f\|_\infty \hat\m B_\e(\partial R_0).
\end{eqnarray}

We may assume that $\hat \m R_0 + \hat \m B_\e(\partial R_0)\le 1$.
By Proposition~\ref{properties2}(i) and (\ref{falphanorm}) we have that
$\e^{m+\a}  \|f\|_\infty \le \c_m^{-1}\e^{\a}\|f\|_\a$
and $\|f\|_\a=|f|_\a+\|f\|_1$ respectively.
So by (\ref{osc(pf)}), (\ref{R1}) and (\ref{R2}), we get
\begin{eqnarray*}
\int_{\hat M} \osc\bigl(Pf, \ B_\e(\cdot)\bigr)d\hat \m
\!\!&\le&\!\!
  \bigl[ (1+\dd \e^\a) \bigl(s^\a\e^\a+ 2G(\e, \e_0)\e_0^\a\bigr)
  +3K'b\c_m^{-1}\e^{\a}\bigr]|f|_\a    \nonumber\\
&+&\! \bigl[2\dd \e^\a+ 2(1+\dd \e^\a) G(\e, \e_0)
  +3K'b\c_m^{-1}\e^{\a}\bigr]\|f\|_1.
\end{eqnarray*}
Now the result follows by the choice of $\eta$ and $D$
in (\ref{xib}) and (\ref{DD}).
\qed

{\it Acknowledgment:}
During the paper was written,
H. Hu was partially supported by the National Science Fundation
grants DMS-9970646 and DMS-0240097.  He also wishes to acknowledge
the hospitality and support of the Centre de Physique Th\'eorique,
CNRS, Marseille, where part of work was done.
S. Vaienti thanks the Department of Mathematics at Michigan state
University for the kind hospitality during the completation of this work.


\end{document}